\documentclass[10pt]{amsart}

\usepackage{amssymb}
\usepackage{amsmath}
\usepackage{amsfonts}
\usepackage{hyperref}
\usepackage{multirow, comment}
\usepackage{mathrsfs}
\usepackage[mathscr]{euscript}

\theoremstyle{plain}
\newtheorem{theorem}{Theorem}[section]
\newtheorem{proposition}[theorem]{Proposition}

\newtheorem{lemma}[theorem]{Lemma}

\theoremstyle{definition}
\newtheorem{remark}[theorem]{Remark}
\newtheorem{definition}[theorem]{Def\/inition}

\newcommand{\begpf}{\noindent{\bf Proof.}\enspace}
\newcommand{\epf}{{\ifhmode\unskip\nobreak\hfil\penalty50 \hskip1em
\else\nobreak\fi \nobreak\mbox{}\hfil\mbox{$\square$} \parfillskip=0pt
\finalhyphendemerits=0 \par\vskip5pt}}
\newcommand{\hra}{\hookrightarrow}
\newcommand{\iso}{\ \smash{\mathop{\longrightarrow}\limits^{\sim}}\ }
\newcommand{\lra}{\longrightarrow}
\newcommand{\ol}{\overline}

\newcommand{\A}{\mathbf{A}}
\newcommand{\C}{\mathbf{C}}
\newcommand{\F}{\mathbf{F}}
\newcommand{\N}{\mathbf{N}}
\newcommand{\Q}{\mathbf{Q}}
\newcommand{\R}{\mathbf{R}}
\newcommand{\T}{\mathbf{T}}
\newcommand{\Z}{\mathbf{Z}}
\DeclareMathOperator{\GL}{GL}

\DeclareMathOperator{\PGL}{PGL}

\newcommand{\Ind}{\mathrm{Ind}}

\newcommand{\Sym}{\mathrm{Sym}}
\renewcommand{\det}{\mathrm{det}}
\newcommand{\res}{\mathrm{res}}
\newcommand{\Res}{\mathrm{Res}}
\newcommand{\sgn}{\mathrm{sign}}
\newcommand{\h}{\mathrm{\mathfrak{h}}}

\newcommand{\gal}{\mathrm{Gal}}
\newcommand{\tr}{\mathrm{tr}}

\newcommand{\lcm}{\mathrm{lcm}}

\newcommand{\FF}{\mathbf{F}}
\newcommand{\CO}{\mathcal{O}}

\newcommand{\ga}{\mathfrak{a}}

\newcommand{\gn}{\mathfrak{n}}
\newcommand{\gm}{\mathfrak{m}}

\newcommand{\Frob}{\mathsf{Frob}}

\newcommand{\Fpbar}{{\overline{\F}_p}}

\begin{document}

\title[Crystalline lifts of automorphic Galois representations]{Crystalline lifts of 
two-dimensional mod $p$ automorphic Galois representations}

\author{\sc Fred Diamond}
\address{Department of Mathematics,
King's College London, WC2R 2LS, UK}
\email{fred.diamond@kcl.ac.uk}
\author{\sc Davide A. Reduzzi}
\address{Department of Mathematics,
The University of Chicago, Chicago, IL 60637, USA}
\email{reduzzi@math.uchicago.edu}

\thanks{F.D.~was partially supported by a Leverhulme Trust Research Project Grant, and partially 
by EPSRC Grant EP/L025302/1.  D.R.~was partially supported by an
AMS-Simons Research Travel Grant.}
\date{September 2015}

\begin{abstract}
We show that a sufficient condition for an irreducible automorphic Galois representation $\rho:
G_F\to\GL_2(\Fpbar)$ of a totally real field $F$ to have an automorphic crystalline lift is that for each place $v$ of $F$ above $p$ the restriction $\det\rho|_{I_v}$ is a fixed power of the mod $p$ cyclotomic character. Moreover, we show that the only obstruction to controlling the level and character of such automorphic lifts arises for badly dihedral representations.
\end{abstract}

\maketitle

\section{Introduction}  \label{sec:intro}
Let $\rho:G_{\Q}\to\GL_2(\ol\F_p)$ be a continuous irreducible odd representation of the absolute Galois group of the rationals. By Serre's conjecture, now a theorem of Khare-Wintenberger and Kisin, $\rho$ is automorphic. It is moreover known that $\rho$ arises from a modular form of level prime to $p$.

The analogue of the last statement for automorphic mod $p$ representations of the absolute Galois group of a totally real number field is false in general. The purpose of this note is to give sufficient conditions for a
mod $p$ Hilbert modular form to have a lift of level prime to $p$, or equivalently
for the associated Galois representation to have an automorphic lift
which is crystalline at all primes over $p$.  The result is motivated
by a question of Dieulefait and Pacetti; see \cite[Lemma~8.32]{dp}
where a special case of Theorem~\ref{main} is used in the construction
of ``chains'' of compatible systems of Galois representations.

We first recall the existence of an obvious obstruction. Let $p$ be a prime number and $F$ a totally real number field. Denote by $\Sigma$ the set of embeddings of $F$ in $\R$. For integers $k\geq 2$ and $w$ having the same parity denote by $D_{k,w}$ the discrete series representation of $\GL_2(\R)$ having Blatter parameters $(k,w)$. In particular $D_{k,w}$ has central character $t \mapsto t^{-w}$. Fix a tuple $\vec{k}=(k_\tau)_{\tau\in\Sigma}\in\Z_{\geq 2}^\Sigma$ and an integer $w$ all sharing the same parity. Let $\pi$ be a cuspidal automorphic representation of $\GL_2(\A_F)$ which is holomorphic of weight $(\vec{k},w)$, \emph{i.e.}, such that $\pi_\tau\simeq D_{k_\tau,w}$ for all $\tau\in\Sigma$. Assume moreover that the level of $\pi$ is coprime with $p$. If $\rho_\pi:G_F \to \GL_2(\ol\Q_p)$ denotes the $p$-adic Galois representation attached to $\pi$, we have $\det\rho_\pi|_{I_v}=\epsilon_v^{w-1}$ for all primes $v$ of $F$ dividing $p$, where $I_v$ is the inertia subgroup of a decomposition group of $G_F$ at $v$, and $\epsilon_v$ is the $p$-adic cylotomic character restricted to $I_v$ (cf. \cite[Corollary 2.11]{bdj} and \cite{ca}).

We show that this condition on $\det\rho_\pi|_{I_v}$ is the only obstruction to the existence of  a crystalline lift of an irreducible automorphic representation $\rho:G_F\to\GL_2(\Fpbar)$.  Moreover we control the conductor and central character of such a lift provided
only that $\rho$ is not badly dihedral (see Definition~\ref{defn:bad}).  More precisely, we prove:

\begin{theorem}\label{main}
Suppose that $\rho:G_F \to \GL_2(\Fpbar)$ is automorphic, irreducible,
and that for some integer $k$, we have $\det\rho|_{I_v} = \overline\epsilon_v^{k-1}$ for all
$v|p$.  Then there exists $n_0$ such that if $n\ge n_0$, there exists a cuspidal
automorphic representation $\pi$ of $\GL_2(\A_F)$ such that 
\begin{itemize}
\item if $v|p$, then $\pi_v$ is unramified principal series;
\item if $v|\infty$, then $\pi_v \cong D_{k+n\delta,k+n\delta}$ where $\delta = \lcm\{\,(p-1)/\gcd(p-1,e_v)\,|\, v|p\,\}$;
\item $\overline{\rho}_\pi \cong \rho$.
\end{itemize}
Suppose further that $\rho$ is automorphic with prime-to-$p$ conductor dividing
$\gn\subset\CO_F$, and that $\psi$
is a finite order Hecke character of $\A_F^\times$ of 
conductor dividing $\gn$, totally of parity  $w = k+n\delta$,
and satisfying $\det\rho = \overline{\psi}\overline{\epsilon}^{w-1}$.
Then if $\rho$ is not badly dihedral, we can choose
$\pi$ as above with conductor dividing $\gn$
and central character $\psi^{-1}|\ |^{-w}$.
\end{theorem}

(Here $e_v$ denotes the ramification degree of $v$ over $p$).
We remark that we have ensured in the conclusion that the lift has parallel
weight since it seems no harder to achieve and slightly simplifies the statement.

There are two main ingredients to the proof of Theorem~\ref{main}.
The first of these is Proposition~\ref{llification} below, which is a statement
purely about mod $p$ representations of $\GL_2(\FF)$ where $\FF$ is a finite
field of characteristic $p$.
The result can be deduced from the main result of \cite{ro}, but we
give instead a short self-contained proof that could be useful if one wishes
to extract an explicit value of $n_0$ in the conclusion of Theorem~\ref{main}.

We will then deduce Theorem~\ref{main} from standard arguments for
producing congruences and liftings of cohomology classes (cf. section \ref{First}).  We must
do some work however to show that the only obstruction to controlling the level and
character is in the case of badly dihedral representations, even for $p=2$ (cf.~sections \ref{Second} and \ref{Second_odd}).   A related result is proved in \cite[Lemma~4.11]{bdj} using the
Galois action on the cohomology of Shimura curves, but since we also wish to work
with forms on definite quaternion algebras, we give a different argument in this paper
by interpreting the obstruction in terms of the Hecke action.
We remark that this obstruction is genuine,
but that even in this case one can obtain slightly weaker results by modifying
our arguments or by constructing CM lifts.
We remark also that we could instead have attempted to deduce a version
of Theorem~\ref{main} using level lowering or automorphy lifting
theorems, as in \cite{ja} or \cite{ge}, but this approach would have required
additional hypotheses, such as adequacy of the image of $\rho$,
and made the case of $p=2$ even more problematic.

Finally we also give two more refined variants of the main result (Theorems~\ref{thm:variant1}
and~\ref{thm:variant2}) in the special case where the initial automorphic representation
has weight $(2,\ldots,2)$ and is unramified or special at all primes over $p$.  This is
again with a view to applications along the lines of those in \cite{dp}.

\subsection*{Acknowledgements}  We thank Luis Dieulefait for calling our attention to
the question addressed in this paper, and Lassina Demb\'el\'e for providing the
example in Remark~\ref{lassina}.  We are also grateful to MSRI and the organizers
of its program on New Geometric Methods in Number Theory and Automorphic Forms,
where some of the work was carried out.

\section{Grothendieck ring relations}  \label{sec:relations}
In this section we denote by $\FF$ a fixed finite field of characteristic $p>0$. Fix an
embedding $\tau_0: \FF \to \Fpbar$ and let $\tau_i = \tau_0\circ \Frob^i$
where $\Frob$ is the (arithmetic) Frobenius automorphism of $\FF$ and we 
sometimes view $i \in \Z/f\Z$ where $f = [\FF:\F_p]$.

For $i = 0,1,\ldots,f-1$ and $n\ge 0$, we let $\Sym^n_{[i]} = 
\Fpbar\otimes_{\FF,\tau_i} \Sym^{n} \FF^2$
where $\Sym^{n}\FF^2$ denotes the $n$th symmetric power of the
standard representation of $\GL_2(\FF)$ for $n\ge 0$; by convention we let
$\Sym_{[i]}^{-1}\FF^2 = 0$.
For $\vec{n} = (n_0,n_1,\ldots,n_{f-1})$ with $n_0,\ldots,n_{f-1} \ge -1$, we let 
$$S_{\vec{n}} = \otimes_{i=0}^{f-1}\Sym_{[i]}^{n_i}.$$
Recall that $S_{\vec{n}}$ is an irreducible representation of $\GL_2(\FF)$ if and only if
$0\le n_i \le p-1$ for all $i$, and that every irreducible representation 
of $\GL_2(\FF)$ over $\Fpbar$ is of the form $\det^a \otimes S_{\vec{n}}$ for
some $a\in \Z$ and $\vec{n}$ as above.

We let
$G_0(\Fpbar[\GL_2(\FF)])$ denote the Grothendieck group on
finite-dimensional representations of $\GL_2(\FF)$ over $\Fpbar$,
which is thus isomorphic to the free abelian group generated
by the classes $[\det^a\otimes S_{\vec{n}}]$ for $a = 0,\ldots,p^f-2$
and $\vec{n} = (n_0,n_1,\ldots,n_{f-1})$ as above.

We use $\leq$ for the natural partial ordering on $G_0(\Fpbar[\GL_2(\FF)])$;
thus $R \le R'$ whenever $R' - R$ is in the submonoid
of $G_0(\Fpbar[\GL_2(\FF)])$ consisting of classes of (actual) 
$\Fpbar$-representations of $\GL_2(\FF)$.  Note that if
$\sigma$ and $\sigma'$ are $\Fpbar$-representations of $\GL_2(\FF)$, then
$[\sigma] \le [\sigma']$ if and only if there is an embedding of the semisimplification
of $\sigma$ in that of $\sigma'$.  In particular if $\sigma$ is irreducible, then 
$[\sigma] \le [\sigma']$ if and only if $\sigma$ is a Jordan-H\"older factor of $\sigma'$.

Assume $f \geq 1$ is arbitrary. If $k<-1$
for each $i\in\Z/f\Z$ define (cf. \cite{se}):
\[
\left[\Sym_{[i]}^k\right]:=-\left[\det^{p^i(k+1)}\otimes\Sym^{-k-2}_{[i]} \right].
\]
In what follows we slightly abuse notation by allowing taking brackets of virtual representations 
in $G_0(\Fpbar[\GL_2(\FF)])$.

Denote by $\N_{\FF/\F_p}$ the field norm map for the extension $\FF/\F_p$.

\begin{proposition} \label{llification}
Let $\sigma$ be an irreducible representation of
$\GL_2(\FF)$ over $\Fpbar$ with central character of the form $\N_{\FF/\F_p}^{es}$
for some $e,s \ge 1$.  Then $[\sigma]\le \left[S_{(t,\dots,t)}^{\otimes e}\right]$
for all sufficiently large $t \equiv s \bmod (p-1)/\gcd(p-1,e)$.
\end{proposition}
The proof will be based on the following lemmas, which can be
viewed as providing algebraic analogues of theta operators 
and Hasse invariants in order to shift weights of automorphic
forms in characteristic $p$.

Let $n,m \geq 0$. The usual identification of graded $\FF$-algebras $\Sym\,\FF^2\simeq \FF[t_1,t_2]$
induces an action of $\GL_2(\FF)$ on $\FF[t_1,t_2]$. When $f=1$, multiplication by the 
Dickson invariant $t_1^p t_2-t_1 t_2^p\in \FF[t_1,t_2]^{\mathrm{SL}_2(\FF)}$ 
induces a $\GL_2(\FF)$-equivariant embedding $S_n \hra \det^{-1} \otimes S_{n+p+1}$ (\cite[section 3]{as}).
Similarly, when $f>1$ we obtain a $\GL_2(\FF)$-equivariant embedding
$\Sym_{[0]}^{m}\otimes \Sym_{[1]}^{n}\hra
\det^{-p}\otimes \Sym_{[0]}^{m+p}\otimes \Sym_{[1]}^{n+1}$ induced by $t_1^p\otimes t_2-t_2^p\otimes t_1$
 (\cite[section~3.5.1]{re}). We obtain in particular:

\begin{lemma} \label{theta:Fp}
Suppose $f=1$ and $n\ge 0$. Then $[S_n] \le [\det^{-1} \otimes S_{n+p+1}]$.
\end{lemma}

\begin{lemma}\label{theta:Fq}
 Suppose $f > 1$ and $m,n \ge 0$. Then
 $$\left[\Sym_{[0]}^{m}\otimes \Sym_{[1]}^{n}\right] \le \\ 
\left[\det^{-p}\otimes  \Sym_{[0]}^{m+p}\otimes \Sym_{[1]}^{n+1}\right].$$
\end{lemma}

\begin{lemma} \label{hasse:Fp}
Suppose $f = 1$ and $n \ge 0$. Then $[S_n] \le [S_{n+p-1}]$ unless $n = r(p+1)$ for some $r$, in which case
we have $[S_n] \le [S_{n+p-1} + \det^r]$.
\end{lemma}

\begpf
Serre's periodic relation: $$[S_{n+p-1}-S_n]=[\det\otimes(S_{n-2}-S_{n-p-1})],$$ valid in $G_0(\Fpbar[\GL_2(\FF)])$ for all $n\in\Z$ (cf. \cite{se})
implies that the positivity of $[S_{n+p-1}-S_n]$ depends only on $n \mod p+1$. If $1\leq
n < p+1$ the term $[S_{n-p-1}]$ is non-positive, so that $[S_n] \leq[S_{n+p-1}] $ when $n\not\equiv 0 \mod p+1$.
When $n=r(p+1)$ we see by induction on $r$ that $[S_{n+p-1}-S_n]= [\det^{r}\otimes
(S_{p-1} - 1)]$.
\epf

\begin{lemma} \label{hasse:Fq}
Suppose $f>1$ and $np > m \ge 0$. Then 
$$\left[\Sym_{[0]}^m\otimes\Sym_{[1]}^n\right] \le \\
\left[\Sym^{m+p}_{[0]}\otimes \\ \Sym_{[1]}^{n-1}\right].$$
\end{lemma}

\begpf
We proceed by induction on $n$. The statement for $n=1$ follows from the identity $\left[\Sym^{m+p}_{[0]}\right]=
\left[\Sym_{[0]}^m\otimes\Sym_{[1]}^1-\det^p\otimes\Sym_{[0]}^{m-p}\right]$ (cf. last equation in \cite[Theorem~2.7]{re}),
together with the fact that $\left[\Sym_{[0]}^{m-p}\right] \leq 0$ since $m<p$. 
Assuming now the statement 
for a fixed $n=n_0>0$ and letting $0\leq m < (n_0+1)p$, we have, using again the above identity:
\begin{equation*}
\begin{split}
\left[\Sym^{m+p}_{[0]}\otimes\Sym_{[1]}^{n_0}\right]=\left[\Sym^m_{[0]}\otimes\Sym^1_{[1]}\otimes\Sym^{n_0}_{[1]}
-\det^p\otimes\Sym_{[0]}^{m-p}\otimes\Sym_{[1]}^{n_0}\right]\\
=\left[\Sym^m_{[0]}\otimes\Sym^{{n_0}+1}_{[1]}+\det^p\otimes
\left(\Sym^m_{[0]}\otimes\Sym^{{n_0}-1}_{[1]}-
\Sym^{m-p}_{[0]}\otimes\Sym^{n_0}_{[1]}\right)\right],
\end{split}
\end{equation*}
If $m<p$, then $\left[\Sym^{m-p}_{[0]}\right]\le 0$ in $G_0(\Fpbar[\GL_2(\FF)])$, so the expression
between rounded parenthesis is positive, implying the statement for $n_0+1$. If $m\geq p$, then 
$0\leq m-p < n_0 p$ so that $
\left[\Sym^{m-p}_{[0]}\otimes\Sym^{n_0}_{[1]}\right]\leq\left[\Sym^m_{[0]}\otimes\Sym^{{n_0}-1}_{[1]}\right]$. 
The result follows.
\epf

We now proceed with the proof of Proposition~\ref{llification}.
Suppose that $\sigma = \det^{a} \otimes S_{\vec{n}}$ with $a\ge 0$ and
$\vec{n}=(n_0,n_1,\dots,n_{f-1})$ with $0\leq n_i\leq p-1$.

We first treat the case $e= f =1$.   By Lemma~\ref{theta:Fp} and induction
on $n$ , we have $[\sigma] \le [S_{n+a(p+1)}]$.   Let $t_0 = n + a(p+1)$, and
note that $t_0 \equiv s \bmod p - 1$.   If $n=0$, then we may replace $a$
by $a + p - 1$ so as to assume $t_0 \ge p^2 - 1$.  We claim that 
$[\sigma] \le [S_t]$ for all $t \equiv s \bmod p - 1$ with $t \ge t_0$.
Indeed it suffices to prove that if $t\ge t_0$ and $[\sigma] \le [S_t]$,
then $[\sigma] \le [S_{t+p-1}]$, and this is immediate from 
Lemma~\ref{hasse:Fp} except in the case $\sigma = \det^r$,
$t = r(p+1)$.   Note however that Lemma~\ref{theta:Fp} implies
that $(b+1)[\det^u] \le [S_{u(p+1)}]$ if $u \ge (p-1)b$; in particular
$2[\det^r] \le [S_t]$ if $t = r(p+1) \ge p^2-1$, so in this case it again
follows from Lemma~\ref{hasse:Fp} that $[\sigma] = [\det^r] \le [S_{t+p-1}]$.

Next we treat the case $e=1$, $f > 1$.   Note that Lemma~\ref{theta:Fq}
implies that  $\left[\Sym_{[i]}^{m}\otimes \Sym_{[i+1]}^{n}\right]\le
\left[\det^{-p^{i+1}}\otimes \Sym_{[i]}^{m+p}\otimes \Sym_{[i+1]}^{n+1}\right],$
for any $m,n\ge 0$, $i\in \Z/f\Z$, and hence that
\begin{equation}\label{bigtheta}
[S_{\vec{m}}] \le \left[\det^{-\sum_{i=0}^{f-1} b_ip^i} \otimes S_{\vec{m}'}\right]
\end{equation}
for any $m_0,\ldots,m_{f-1},b_0,\ldots,b_{f-1} \ge 0$, where
$\vec{m} = (m_0,\ldots,m_{f-1})$ and $\vec{m}' = (m_0 + b_0 + pb_1,\ldots,m_{f-1} +b_{f-1}
 + pb_0)$.  In particular writing $a = \sum_{i=0}^{f-1}b_ip^i$ with $b_0,\ldots,b_{f-1} \ge 0$,
 we have that  $[\sigma] \le [S_{\vec{n}'}]$ where 
 $\vec{n}' = (n_0',\ldots,n_{f-1}')$ with $n'_i =  n_i + b_i + pb_{i+1}$.
 Note that
 $$\sum_{i=0}^{f-1} n_i'p^i \equiv 2a + \sum_{i=0}^{f-1} n_ip^i \equiv 
       s\left(\sum_{i=0}^{f-1}p^i \right) \bmod p^f - 1,$$
 from which it follows that $\sum_{i=0}^{f-1} n'_{i+j} p^i$ is divisible
 by $\sum_{i=0}^{f-1}p^i$ for $j=0,\ldots,f-1$.

Now consider the system of equations
\begin{equation}\label{system}
n_0' - x_0 +px_1 = n_1' - x_1 + px_2 = \cdots = n'_{f-1} - x_{f-1} + px_0.
\end{equation}
For any $x_0 \in \Z$, we obtain a solution with $x_0,\ldots,x_{f-1} \in \Z$ by
setting
$$ x_j = x_{j-1}  - n'_{j-1} + 
\left.\left(\sum_{i=0}^{f-1} n'_{i+j} p^i \right)\right/\left(\sum_{i=0}^{f-1}p^i\right)$$
for $j=1,\ldots,f-1$.
In particular we may choose a solution of (\ref{system})
with $x_0,\ldots,x_{f-1}$ non-negative integers.

We now wish to apply Lemma~\ref{hasse:Fq}, or rather its twist by $\Frob^i$,
iteratively $x_{i+1}$ times for $i = 0,\ldots,f-1$ in order to conclude that
 $[\sigma] \le [S_{(t,\ldots,t)}]$ where $t$ is the common value of
$n_i' - x_i + px_{i+1}$, but we must first ensure that the inequality $np > m$
in the hypothesis of the lemma is satisfied at each stage.  To this
end, note that we may replace $\vec{n}'$ by
$\vec{n}'' = \vec{n}' + r(p^2-1,\ldots,p^2-1)$ for any integer $r\ge 0$;
indeed  $[\sigma] \le [S_{\vec{n}'}] \le [S_{\vec{n}''}]$ by
(\ref{bigtheta}), and (\ref{system}) still holds with each $n'_i$ replaced
by $n''_i = n'_i + r(p^2-1)$.   Choosing $r$ so that 
$$r(p-1)(p^2-1) > n_i' - pn_{i+1}' +2px_{i+1}$$
for each $i$, we find that $p(n''_{i+1} - x_{i+1}) > n''_i + px_{i+1}$.
Now by Lemma~\ref{hasse:Fq} and induction on $\sum_{i=0}^{f-1}d_i$,
we see that if $0 \le d_i \le x_i$ for $i = 0,\ldots,f-1$, then
$[S_{\vec{n}''}] \le [S_{\vec{n}'''}]$ where $n_i''' = n_i'' - d_i + pd_{i+1}$.
It follows that $[\sigma] \le [S_{(t,\ldots,t)}]$ where $t$ is the common
value of $n_i'' - x_i + px_{i+1}$.  Similarly we find that if $t > 0$,
then $[S_{(t,\ldots,t)}] \le [S_{(t+p-1,\ldots,t+p-1)}]$, which completes
the proof in the case $e=1$, $f>1$.

Finally we treat the case $e > 1$.  From the case $e=1$, we have that
$[\sigma] \le [S_{(u,\ldots,u)}]$ for all sufficiently large $u \equiv es \bmod p-1$,
hence $[\sigma] \le [S_{(et,\ldots,et)}]$ for all sufficiently large 
$t\equiv s \bmod (p-1)/\gcd(e,p-1)$.   From the natural surjection
$$\left( \Sym_{[i]}^t \right)^{\otimes e}  \to \Sym_{[i]}^{et}$$
we see that $[S_{(et,\ldots,et)}] \le \left[S_{(t,\ldots,t)}^{\otimes e}\right]$,
concluding the proof of Proposition~\ref{llification}. \\

\begin{remark}  
When $f=1$, multiplication by the Dickson invariant
$(t_1^{p^2}t_2-t_1 t_2^{p^2})/(t_1^p t_2-t_1 t_2^p) \in (\Sym\,\FF^2)^{\GL_2(\FF)}$ induces
a $\GL_2(\FF)$-equivariant injection $S_k\to S_{k+p(p-1)}$ for all $k\ge 0$. Notice that the change in weight produced 
by this operator does not allow us to prove the desired result.
\end{remark}

\section{Lifting to characteristic zero}  \label{sec:lifting}
The first part of Theorem \ref{main} is proved in section \ref{First} below. In sections \ref{Second} and \ref{Second_odd} we refine the argument to control the level and character of the crystalline lifts we produce, thus proving the second part of Theorem \ref{main}. We begin by fixing some notation.

\subsection{Notation}
We normalize local and global class field theory so that geometric Frobenius elements correspond to uniformizers, and we adopt Hecke's normalizations of local-global compatibility when associating a Galois representation to an automorphic form. These are the normalizations adopted in \cite{ca} and \cite{bdj}. 

Let $p$ be a prime number. We fix an algebraic closure $\ol\Q$ (resp. $\ol\Q_p$) of the field $\Q$ of rational numbers (resp. of the field $\Q_p$ of $p$-adic numbers). We choose an embedding $\ol\Q\to\C$ and an isomorphism $\ol\Q_p\cong\C$, so that we also identify $\ol\Q$ with a subfield of $\ol\Q_p$. Denote by $\ol\FF_p$ a fixed algebraic closure of the field $\FF_p$ with $p$ elements.

Let $F\subset\ol\Q$ be a totally real number field. Denote by $G_F=\mathrm{Gal}(\ol\Q/F)$ its absolute Galois group, and by $\epsilon:G_F\to\Z_p^\times$ the $p$-adic cyclotomic character. Let $\Sigma$ be the set of embeddings of $F$ in $\ol\Q$. 

For each finite place $x$ of $F$ we denote by $F_x$ the completion of $F$ at $x$, and by $\CO_{F_x}$ its ring of integers. We let $\A_F$ (resp. $\A_{F,f}$) denote the topological ring of ad\`eles (resp. finite ad\`eles) of $F$.

Let $v$ be a place of $F$ lying above $p$. Let $G_{v}$ denote the decomposition group of $G_F$ at $v$ induced by $\ol\Q\subset\ol\Q_p$, and $I_v$ be its inertia subgroup. We let $k_v$ be the residue field of $\CO_{F_v}$, and we set $f_v:=[k_v:\F_p]$. Denote by $e_v$ the absolute inertial degree of $v$, and by $\overline{\Sigma}_v$ be the set of embeddings of $k_v$ in $\Fpbar$. We let $\epsilon_v$ denote the restriction of the $p$-adic cyclotomic character to $G_{v}$ or to $I_v$. The reduction modulo $p$ of $\epsilon_v$ is denoted by $\ol\epsilon_v$. 

For any $\tau\in\overline{\Sigma}_v$ we denote by $\omega_\tau$ the corresponding fundamental
character of $I_v$, defined as the composition 
$I_v\lra\CO_{v}^\times\lra k_v^\times \overset{\tau}\lra \Fpbar^\times$, where the first map is the restriction of the inverse of the reciprocity isomorphism of local class field theory. Recall that the restriction to $I_v$ of the local mod $p$ cyclotomic character $\overline\epsilon_v$ of $G_{v}$ is
given by $\prod_{\tau\in\overline{\Sigma}_v}\omega_\tau^{e_v}$.

For integers $k\geq 2$ and $w$ having the same parity we denote by $D_{k,w}$ the discrete series representation of $\GL_2(\R)$ having Blatter parameters $(k,w)$, and hence central character $t\mapsto t^{-w}$.

\subsection{Existence of lifts}\label{First}
For each place $v|p$ of $F$, fix an embedding $\tau_{v,0}:
k_v\to\Fpbar$ and, as in Section \ref{sec:relations}, set $\tau_{v,i}=\tau_{v,0}\circ\Frob^i_v$ where $\Frob_v$ denotes the arithmetic
Frobenius of $k_v$ and $i\in\Z /f_v\Z$. 
For any $\vec{n}=(n_0,\dots,n_{f_v-1})$ with $n_0,\dots,n_{f_v-1}\geq -1$ 
let 
$$S_{v,\vec{n}}=\bigotimes_{i=0}^{i=f_v-1} \left(\Fpbar\otimes_{k_v,\tau_{v,i}}\Sym^{n_i}k_v^2\right),$$
viewed as an $\Fpbar$-linear representation of $\GL_2(k_v)$.

Suppose now that $\rho:G_F\to\GL_2(\Fpbar)$ is the modular Galois representation 
in the statement of Theorem \ref{main}. Assume that $\sigma$ is a Serre weight for $\rho$ in the
sense of \cite{bdj}; in particular $\sigma$ is an $\Fpbar$-linear irreducible representation
of $\GL_2(\CO_F/(p))=\prod_{v|p}\GL_2(\CO_F/(v^{e_v}))$, and therefore it can be written
as $\sigma=\otimes_{v|p}\sigma_v$ where each $\sigma_v$ is an irreducible representation
of $\GL_2(k_v)$. 

Denote by $\N_{k_v/\F_p}$ the norm map attached to the field extension $k_v/\F_p$. Slightly modifying the proof of \cite[Corollary 2.11]{bdj} by taking the map $N:(\CO_F/(p))^\times\to\F_p^\times$
considered there to be $\prod_{v|p}\N_{k_v/\F_p}^{e_v}$, and by applying the generalization 
to the ramified settings of
\cite[Proposition 2.10]{bdj}, we deduce that if the central character of $\sigma_v$
is given by $\prod_{\tau\in\overline\Sigma_v}\tau^{c_\tau-2e_v }$ for some integers $c_\tau$, then
\[
\det\rho|_{I_v}=\prod_{\tau\in\overline\Sigma_v}\omega_\tau^{c_\tau-e_v}.
\]

By assumption, for any prime $v$ of $F$ above $p$ we have 
$\det\rho|_{I_v}=\overline\epsilon^{k-1}_v=\prod_{\tau\in\overline{\Sigma}_v} \omega_\tau^{e_v (k-1)}$. We deduce therefore that
 the central character of $\sigma_v$
is given by:
\[
\N_{k_v/\F_p}^{e_v (k-2)}.
\]
Proposition \ref{llification} implies that $\sigma_v$ is a 
Jordan-H\"older factor of $S_{v,(t_v,\dots,t_v)}^{\otimes e_v}$ for all sufficiently large $t_v\equiv k-2 \bmod
(p-1)/\gcd(p-1,e_v)$. Define $\delta := \lcm\{\,(p-1)/\gcd(p-1,e_v)\,|\, v|p\,\}$. We can thus find a non-negative integer $n_0$ such that for all $n\geq n_0$ we have $k-2+n\delta \geq 0$ and the weight $\sigma$ 
is a Jordan-H\"older factor of the $\GL_2(\CO_F/(p))$-representation 
$\otimes_{v|p}S_{v,(k-2+n\delta,\dots,k-2+n\delta)}^{\otimes e_v}$. By (a generalization 
to the ramified settings of) \cite[Proposition 2.5]{bdj} we deduce
that $\rho$ arises from a cuspidal automorphic representation $\pi$ of $\GL_2(\A_F)$ having level prime to $p$ and 
such that $\pi_v\simeq D_{k+n\delta,k+n\delta}$ for all places $v|\infty$ of $F$. This proves the first part of Theorem \ref{main}.\\

\subsection{The refinement of the argument in the case $[F:\Q]$ even}\label{Second}

We now explain how to refine the argument above in order to control the level and
character of $\pi$, and prove the second half of Theorem \ref{main}.   The case of $[F:\Q]$ odd can be treated by modifying the argument
above and using results in the proof of \cite[Lemma~4.11]{bdj} to show that
obstructions arise only for badly dihedral representations.    For the case
of $[F:\Q]$ even, we will need to use forms on definite quaternion algebras and
prove analogous results concerning the obstructions, which we proceed to do
first.

\subsubsection{Automorphic forms on definite quaternion algebras}\label{Sub:DQA}
Suppose that $[F:\Q]$ is even and
let $D$ be the totally definite quaternion algebra over $F$ which splits at all finite places of $F$. Let $\CO_D$ be a fixed maximal order in $D$, and choose isomorphisms of $\CO_{F_x}$-algebras $\CO_{D,x}\cong M_2(\CO_{F_x})$ for each finite place $x$ of $F$. Let $U=\prod_xU_x$ be an open
compact subgroup of $D_f^\times:=(D\otimes_F\A_{F,f})^\times$ such that $U_x\subset \CO_{D,x}^\times$ for each finite place $x$. Let $A$ denote the field $\Fpbar$ or a topological $\Z_p$-algebra of finite type, and fix a continuous representation of  $U\A_{F,f}^\times/F^\times$ 
on a finitely generated (topological) $A$-module $V$.
Let 
$$S_V(U) = \{\ f: D_f^\times \to V \,|\,  f(\gamma g u) = u^{-1} f(g) \mbox{ for
     all $\gamma \in D^\times$, $g \in D_f^\times$, $u \in U\A_{F,f}^\times$}\,\}.$$
Write $D_f^\times=\coprod_{i\in I}D^\times t_i U\A_{F,f}^\times$ where $I$ is a finite set, and let
$\Gamma_i$ denote the finite group $F^\times\backslash(U\A_{F,f}^\times \cap t_i^{-1} D^\times t_i)$,
so that we have an isomorphism of $A$-modules:
\begin{equation}\label{H^0}
S_V(U) \overset{\simeq}\lra \oplus_{i\in I} V^{\Gamma_i}
\end{equation}
induced by $f \mapsto \oplus_{i\in I} f(t_i)$.

Let $S$ be a finite set of finite places of $F$ 
containing the places dividing $p$ and the places $x$ such that $U_x$ is not maximal.  
Let $U_S:=\prod_{x\in S}U_x$ and suppose further that the action of $U$ on $V$ factors
through the projection to $U_S$.  
For $x\not\in S$ fix a choice of uniformizer $\varpi_x$ of $\CO_{F_x}$, and write $U_x\Pi_x U_x=\coprod_\alpha h_\alpha U_x$,
where $\Pi_x=\left(\begin{smallmatrix} \varpi_x&0\\ 0&1 \end{smallmatrix}\right)\in \GL_2({F_x})\cong (D\otimes_F{F_x})^\times$. We define the
Hecke operator $T_x$ acting on $f\in S_V(U)$ by 
$$(T_xf)(g):=
\sum_\alpha f(gh_\alpha)$$ 
for all $g\in D_f^\times$. The Hecke algebra ${\bf T}^{{S}}_A:=A[T_x:x\not\in S]$ acts on $S_V(U)$. With a slight abuse of notation, we will often not indicate the weight and level of automorphic forms on which ${\bf T}^{{S}}_A$ acts.

Let $\vec{k}=(k_\tau)_{\tau\in\Sigma}\in\Z_{\geq 2}^\Sigma$ and $w \in \Z$ 
such that $k_\tau \equiv w \bmod 2$ for all $\tau \in \Sigma$, and let
$\psi$ a finite order Hecke character of $\A_F^\times$,
totally of parity $w$, so $\psi(x) = \psi_f(x_f) \prod_{\tau|\infty}\sgn(x_\tau)^w$ for some
character $\psi_f$ of $\A_{F,f}^\times$.  Let $E \subset \ol{\Q}_p \cong \C$ be a sufficiently
large finite extension of $\Q_p$; in particular we assume $E$ contains the values of $\psi$ and the images of all the embeddings $F\to\ol{\Q}_p$. We suppose further that
for each embedding $\tau: F \to E$, we have a splitting $D \otimes_{F,\tau} E \cong M_2(E)$
so that if $v$ is the place of $F$ induced by $\tau$, then the projection of $U$ 
to $(D\otimes_F{F_v})^\times$
is contained in $\GL_2(\CO_E)$.   Thus for each $\tau\in \Sigma$, we obtain a
map $U \to \GL_2(\CO_E)$, and hence an action of $U$ on 
$\det^{(w-k_\tau+2)/2}\otimes\Sym^{k_\tau-2} \CO_E^2$. 
Suppose now that  $\psi_f$ is trivial on
$U \cap \A_{F,f}^\times$, and let $V_{\vec{k},w,\psi}$ be the 
representation of $U\A_{F,f}^\times$ whose restriction to $U$ is defined by
$$\otimes_{\tau\in \Sigma}  (\det^{(w-k_\tau+2)/2}\otimes\Sym^{k_\tau-2} \CO_E^2),$$ 
and whose restriction to $\A_{F,f}^\times$ is defined by the character
$x \mapsto \N(x_p)^w |x|^w \psi_f(x)$. 

We write $S_{\vec{k},w,\psi}(U)$ for
$S_{V_{\vec{k},w,\psi}}(U)$, and we define $S_{\vec{k},w,\psi}^{\mathrm{triv}}(U):=\{0\}$, unless $\vec{k}=\vec{2}$,
in which case we let $S_{\vec{k},w,\psi}^{\mathrm{triv}}(U)$ consist of those functions in $S_{\vec{k},w,\psi}(U)$ that factor through the reduced norm map $D_f^\times\cong\GL_2(\A_{F,f})\overset{\det}\lra\A_{F,f}^\times$.
Setting $S_0 = S_{\vec{k},w,\psi}(U)/S_{\vec{k},w,\psi}^{\mathrm{triv}}(U)$,
we have by the Jacquet-Langlands correspondence that 
$S_0 \otimes_{\CO_{E}} \C \cong \oplus_\pi \pi_f^U$,  the direct sum running over
all holomorphic cuspidal automorphic representations $\pi = \pi_\infty\otimes \pi_f$
of $\GL_2(\A_F)$ such that $\pi_\tau \cong D_{k_\tau,w}$ for all $\tau\in\Sigma$,
and $\pi$ has central character $\psi^{-1}|\ |^{-w}$.

\subsubsection{Conclusion of the argument} \label{conclusion}
We fix an irreducible representation $\rho:G_F\to\GL_2(\ol\F_p)$ as in Theorem \ref{main}. Assume that $\rho$ arises from a holomorphic cuspidal automorphic form $\pi'$ for $\GL_2(\A_F)$ of paritious weight $(\vec{k},w)\in\Z_{\ge 2}^\Sigma\times\Z$, central character $\psi^{-1}|\ |^{-w}$, and level $U=U_pU^p$, where $U_p=\prod_{v| p}U_v$, $U^p=\prod_{x\nmid p}U_x$, and $U_x\subset\GL_2(\CO_{F_x})$ for all $x$.  Let $S$ be a finite set of finite places of $F$ containing the places of $F$ above $p$ and the places at which $\pi'$ is ramified. Let  $\mathfrak{m_\rho}$ denote the 
maximal ideal of the Hecke algebra ${\bf T}^{S}_{\CO_E}$ attached to $\rho$; thus with our
normalizations, $\gm_\rho$ is the kernel of the homomorphism  ${\bf T}^{S}_{\CO_E} \to \ol{\F}_p$
defined by
$$T_x \mapsto \det(\rho(\Frob_x))^{-1}\N(x)\tr(\rho(\Frob_x)) 
= \ol{\psi}(\varpi_x)^{-1}\N(x)^w\tr(\rho(\Frob_x))$$
for all $x \not\in S$.
By what was recalled in \ref{Sub:DQA} we know that
 $S_{\vec{k},w,\psi}(U)_{\mathfrak{m}_\rho}\neq 0$.

In the next section we will prove the following:

\begin{lemma}\label{lem:exactness}
If $\rho$ is not badly dihedral in the sense of Definition \ref{defn:bad} below, then
the functor $\ol V\mapsto S_{\ol V}(U)_{\mathfrak{m}_\rho}$ from finite dimensional $\Fpbar$-vector spaces endowed with a continuous action of $U_S\A^\times_{F,f}/F^\times$ to ${\bf T}^{S}_{\Fpbar,\mathfrak{m_\rho}}$-modules is exact. 
\end{lemma}

 Let $U_\ast:=\GL_2(\CO_F\otimes_\Z\Z_p)
\cdot U^p$, and notice there is a Hecke equivariant injection $S_{\vec{k},w,\psi}(U) \to S_{V'}(U_\ast)$ where $V' = \Ind_{U\A_{F,f}^\times}^{U_\ast\A_{F,f}^\times} V_{\vec{k},w,\psi}$. In particular, we also have that 
$S_{V'}(U_\ast)_{\mathfrak{m}_\rho}\neq 0$. Fix an embedding of the residue field of $\CO_E$ into $\Fpbar$. Using Lemma \ref{lem:exactness} we see that $S_{\sigma}(U_\ast)_{\mathfrak{m}_\rho}\neq 0$
for some Jordan-H\"older constituent $\sigma$ of $V'\otimes_{\CO_E}\Fpbar$ for the action of $\GL_2(\CO_F\otimes_\Z\Z_p)$. 
 The assumption on $\rho|_{I_v}$ implies by Proposition \ref{llification} that $\sigma$ is a constituent of the $\Fpbar$-linear representation 
 $$\otimes_{v|p}  (\det\otimes S_{v,(k'-2,\ldots,k'-2)}^{\otimes e_v})$$
for all sufficiently large
$k' \equiv k \bmod \delta$, where
$\delta$ as in the statement of Theorem~\ref{main}. It follows that $\gm_\rho$ is in the support of $S_W(U_\ast)$, where $W$ is the $\CO_E$-linear representation $V_{\vec{k'},k',\psi}$ (this again uses Lemma \ref{lem:exactness}). Now the result follows by applying the Jacquet-Langlands correspondence to $S_W(U_\ast)\otimes_{\CO_E}\C$ to produce a holomorphic Hilbert modular form with desired weight, level, and central character.

\subsection{Proof of Lemma \ref{lem:exactness}}

We keep the assumptions and notation from the previous section. Suppose that $0\to V\to V_1 \to V_2 \to 0$ is an exact sequence of finite dimensional $\Fpbar$-vector spaces endowed with a continuous action of $U\A_{F,f}^\times$ factoring through $U_S\A_{F,f}^\times$. By \eqref{H^0} we obtain the exact sequence:

\begin{equation}\label{exact}
0\to S_{V}(U)\to S_{V_1}(U) \to S_{V_2}(U) \to \oplus_{i\in I} H^1(\Gamma_i,V),
\end{equation}
where $\Gamma_i=F^\times\backslash(U\A_{F,f}^\times \cap t_i^{-1} D^\times t_i)$ and $D_f^\times=\coprod_{i\in I}D^\times t_i U\A_{F,f}^\times$.
Notice the last term in \eqref{exact} vanishes if $[F(\mu_p):F]>2$, which occurs for example when $p>3$ is unramified in $F/\Q$. We will show that, in general, it vanishes after localization at $\gm_\rho$ if $\rho$ is not badly dihedral.

Note that if we choose another representative $t_i' = \delta t_i u$ for the double coset 
$D^\times t_i U\A_{F,f}^\times$ with $\delta \in D^\times$, $u \in U\A_{F,f}^\times$ and set
$\Gamma_i' = F^\times\backslash(U\A_{F,f}^\times \cap (t_i')^{-1} D^\times t_i')$,
then $\Gamma_i = u\Gamma'_i u^{-1}$ and we obtain canonical isomorphisms
$H^j (\Gamma_i,V) \to H^j(\Gamma_i',V)$ induced by the isomorphisms
$$\begin{array}{ccccccc} 
\Gamma_i'& \to& \Gamma_i, &\quad& V &\to& \Res_{\Gamma_i}^{\Gamma_i'}V \\
                                    g &\mapsto&  ugu^{-1}, & \quad&   v& \mapsto   &uv.\end{array}$$
We let $X_U = D^\times \backslash D_f^\times / U \A_{F,f}^\times$ and
write $H^j(X_U,V)$ for $\oplus_{i\in I} H^j(\Gamma_i,V)$; this is independent of the
choice of the $t_i$ up to canonical isomorphism\footnote{Alternatively 
one can arrive at this notation by defining
a Grothendieck topology on the groupoid fibered over $X_U$ by the $\Gamma_i$ 
and viewing $V$ as a sheaf on the associated site.},
which is moreover compatible with
the isomorphism $S_V(U) \cong H^0(X_U,V)$ in the evident sense.  We may thus
rewrite the exact sequence (\ref{exact}) as
$$0\to H^0(X_U,V) \to H^0(X_U,V_1) \to H^0(X_U,V_2)  \to H^1(X_U,V).$$

\subsubsection{The Hecke action on $H^1(X_U,V)$}\label{HeckeH1}
For $x\notin S$, the Hecke operator $T_x$ acting on $S_V(U)$
can be defined as a composite $\tr \circ (\Pi_x)_* \circ \res$
where $\res$ (resp.~$\tr$) is a restriction (resp.~trace) map to (resp.~from)
forms with respect to a smaller open compact subgroup, and $(\Pi_x)_*$ is
induced by $\Pi_x$. More precisely, consider the natural projection 
$X_{U'} \to X_U$  where $U' = U \cap \Pi_x^{-1} U \Pi_x$, and for each double
coset $D^\times t_i U\A_{F,f}^\times$ in $X_U$, let $\{t_{ij}\}$ be representatives of the preimage
in $X_{U'}$ (so $D^\times t_i U\A_{F,f}^\times$ is the disjoint union over $j$ of the 
$D^\times t_{ij} U'\A_{F,f}^\times$).
Let $\Gamma'_{ij}$ be the corresponding stabilizers, so
$\Gamma'_{ij} = F^\times\backslash(U'\A_{F,f}^\times \cap t_{ij}^{-1}D^\times t_{ij})$.
Writing $t_i  = \delta t_{ij} u$ for some $\delta \in D^\times$, $u \in U'\A^\times_{F,f}$,
we see that $v \mapsto uv$ defines a map $V \to V$ compatible
with the inclusion $\Gamma'_{ij} \to \Gamma_i$ defined
by conjugation by $u$, and this gives a map
$$H^1(\Gamma_i,V) \to \oplus_j H^1(\Gamma'_{ij},V).$$
Taking the direct sum over $i$ of these maps, the resulting map
$\res: H^1(X_U,V) \to H^1(X_{U'},V)$ is independent of the choices of
double coset representatives.

Similarly we define $(\Pi_x)_*$ using the bijection 
$X_{U''} \to X_{U'}$
induced by right multiplication by $\Pi_x$, where $U'' = \Pi_x U'\Pi_x^{-1}$.
Since $D_f^\times = \coprod_{i,j} D^\times t_{ij} U'\A_{F,f}^\times$,
we have $D_f^\times = \coprod_{i,j} D^\times t_{ij}\Pi_x^{-1} U''\A_{F,f}^\times$
and the corresponding stabilizer $\Gamma_{ij}''$ equals $\Pi_x\Gamma_{ij}' \Pi_x^{-1}$.
Since $U_x$ acts trivially on $V$, the isomorphism between the groups
$\Gamma_{ij}'$ and $\Gamma_{ij}''$ defined by conjugation
by $\Pi_x$ is compatible with their action on  $V$,
so it induces an isomorphism $H^1(\Gamma_{ij}',V) \to H^1(\Gamma_{ij}'',V)$.
Taking the direct sum of these isomorphisms gives a well-defined map
$(\Pi_x)_*:  H^1(X_{U'},V) \to H^1(X_{U''},V)$.

Finally $\tr$ is defined similarly to $\res$ but using $X_{U''} \to X_U$
and transfer maps on cohomology.
It is then easy to see that the composite $\tr \circ (\Pi_x)_* \circ\res$ is compatible
with \eqref{exact} and the Hecke operators on the $S_V(U)$ since
each of $\res$, $\Pi_x^*$ and $\tr$ is compatible with \eqref{exact} in the obvious
sense.  Note that in fact $T_x = \tr \circ (\Pi_x)_* \circ\res$ on
$H^0(X_U,V) = S_V(U)$.

More generally if $x_1,x_2,\ldots,x_m$ are distinct primes of $F$ not in $S$,
then we define a Hecke operator $T_{x_1x_2\cdots x_m}$ exactly as
above, but replacing $\Pi_x$ by the product of the $\Pi_{x_i}$, which
we denote by $\Pi_{x_1\cdots x_m}$.
\begin{lemma} \label{lem:product}  We have $T_{x_1x_2\cdots x_m}
 = T_{x_1}T_{x_2}\cdots T_{x_m}$.  In particular the operators $T_{x_i}$
 commute and $\T_{\Fpbar}^S$ acts on $H^1(X_U,V)$.
 \end{lemma}
 \begpf   Consider the diagram:
 $$\begin{array}{ccccccc}
 H^1(X_U,V) & \rightarrow & H^1(X_{U_1'},V) & \rightarrow &
  H^1(X_{U_1''},V) & \rightarrow & H^1(X_U,V) \\&&&&&& \\
 & \searrow & \downarrow & &
 \downarrow && \downarrow \\ &&&&&& \\
 && H^1(X_{U_1' \cap U_2'},V)& \rightarrow &
 H^1(X_{U_1'' \cap U_2'},V) & \rightarrow & H^1(X_{U_2'},V) \\ &&&&&& \\
 &&& \searrow &
 \downarrow && \downarrow \\ &&&&&& \\
 &&&& H^1(X_{U_1'' \cap U_2''}, V) & \rightarrow  & H^1(X_{U_2''},V) \\ &&&&&& \\
 &&&&& \searrow & \downarrow \\ &&&&&& \\
 &&&&&& H^1(X_U,V), \end{array}$$
 where 
 \begin{itemize}
 \item $U_1' = U \cap \Pi_{x_1}^{-1} U \Pi_{x_1}$, $U_1'' = \Pi_{x_1}U_1' \Pi_{x_1}^{-1}$,
 $U_2' = U \cap \Pi_{x_2\cdots x_m}^{-1} U \Pi_{x_2\cdots x_m}$ and
 $U_2'' = \Pi_{x_2\cdots x_m} U_2' \Pi_{x_2\cdots x_m}^{-1}$;
 \item the first row of vertical arrows and the first column of arrows (including the first diagonal) are defined by
  the evident restriction maps;
 \item the last row of arrows and the last column of horizontal arrows (including the last diagonal) are defined by
 the evident  transfer maps;
 \item the middle column (resp. row) of horizontal (resp.~ vertical) arrows is of the form $(\Pi_{x_1})_*$
 (resp.~$(\Pi_{x_2\cdots x_m})_*$) and the middle diagonal arrow is $(\Pi_{x_1x_2\cdots x_m})_*$.
 \end{itemize}
 Note that the top row comprises $T_{x_1}$, the last column $T_{x_2\cdots x_m}$ and the diagonal
 $T_{x_1x_2\cdots x_m}$, so the lemma follows by induction from 
 the commutativity of all the triangles and squares in the diagram.
 
 We only sketch the proof of commutativity of the top right corner, the rest
 being immediate from the definitions of the maps.  Moreover to prove commutativity
 of the top right corner reduces to checking it for the corresponding diagram associated
 to each summand of $H^1(X_{U_2'},V)$.  More precisely, given a double
 coset $D^\times t U_2'\A_{F,f}^\times$ in $X_{U_2'}$ mapping to 
 $D^\times s U \A_{F,f}^\times$ in $X_U$, let $\Delta
  = F^\times\backslash(   U_2'\A_{F,f}^\times\cap t^{-1} D^\times t)$ and
  $\Gamma = F^\times\backslash( U\A_{F,f}^\times \cap s^{-1} D^\times s)$.
  Writing 
  $$D^\times t U_2' \A_{F,f}^\times = \coprod_{j\in J}
  D^\times t_j (U_2'\cap U_1'')\A_{F,f}^\times \quad\mbox{and}\quad
  D^\times s U\A_{F,f}^\times = \coprod_{i\in I}
  D^\times s_i U_1'' \A_{F,f}^\times,$$
  we must check the commutativity of the diagram
  $$\begin{array}{ccc}  \displaystyle{\bigoplus_{i\in I}} H^1(\Gamma_i, V) & \rightarrow & H^1(\Gamma,V)\\
\downarrow && \downarrow \\
  \displaystyle{\bigoplus^{\ }_{j\in J}}H^1(\Delta_j,V) & \rightarrow & H^1(\Delta,V)\end{array}$$
 where $\Delta_j
  = F^\times\backslash(   (U_2'\cap U_1'')\A_{F,f}^\times\cap t_j^{-1} D^\times t_j)$,
  $\Gamma_i
  = F^\times\backslash(  U_1''\A_{F,f}^\times\cap s_i^{-1} D^\times s_i)$
  and the maps are defined as follows:
  \begin{itemize}
  \item Writing $s = \alpha t w$ with $\alpha \in D^\times$, $w \in U\A_{F,f}^\times$,
  the right-hand arrow is the composite
  $$H^1(\Gamma,V) \to H^1(\Delta, \Res_\Gamma^\Delta V) \iso H^1(\Delta,V)$$
  where the inclusion $\Delta \to \Gamma$ is defined by $g \mapsto w^{-1}gw$
  and $\Res_\Gamma^\Delta V \to V$ is defined by $v\mapsto wv$.
  \item  The left-hand arrow is similarly defined component-wise as the composite
 $$H^1(\Gamma_i ,V) \to \bigoplus_j H^1(\Delta_j, \Res_{\Gamma_i}^{\Delta_j} V) 
 \iso \bigoplus_j H^1(\Delta_j,V),$$
 where the direct sum is over $j$ such that $s_i = \alpha_j t_j w_j$ for some
  $\alpha_j \in D^\times$, $w_j \in U_1''\A_{F,f}^\times$.
\item Writing $s = \beta_i s_i y_i$ with $\beta_i \in D^\times$, $y_i \in U\A_{F,f}^\times$
for each $i\in I$,  the top arrow is defined component-wise as the composite
$$H^1(\Gamma_i,V) \iso H^1(\Gamma, \Ind_{\Gamma_i}^\Gamma V)
 \iso H^1(\Gamma,\Ind_{\Gamma_i}^\Gamma\Res_\Gamma^{\Gamma_i} V)
 \to H^1(\Gamma,V),$$
 where the inclusion $\Gamma_i \to \Gamma$ is $g\mapsto y_i^{-1}gy_i$,
 the first isomorphism is that of Shapiro's Lemma, the second
 is induced by $V \iso \Res_\Gamma^{\Gamma_i} V $ defined
 by $v \mapsto y_i^{-1} v$, and the last map is given by the trace
 $\Ind_{\Gamma_i}^\Gamma\Res_\Gamma^{\Gamma_i} V \to V$.
 \item Writing $t = \gamma_j t_j z_j$ with $\gamma_j \in D^\times$,
 $z_j \in U_2'\A_{F,f}^\times$ for each $j$,  the bottom arrow is defined
 similarly by the composite
 $$H^1(\Delta_j,V) \iso H^1(\Delta, \Ind_{\Delta_j}^\Delta V)
 \iso H^1(\Delta,\Ind_{\Delta_j}^\Delta\Res_\Delta^{\Delta_j} V)
 \to H^1(\Delta,V).$$
  \end{itemize}
Note that for each $j \in J$, the resulting diagram of inclusions
$$\begin{array}{ccc}  \Delta_j &\to& \Delta \\ \downarrow&&
\downarrow \\ \Gamma_i & \to & \Gamma \end{array}$$
commutes up to conjugation by the element $g_j = y_i^{-1}w_j^{-1}z_j wF^\times \in \Gamma$.
Unravelling definitions and applying standard functorialities, one is reduced
to checking commutativity of the following diagram of homomorphisms of 
$\Delta$-modules
$$\begin{array}{ccccc}
\displaystyle{\bigoplus_{i\in I}} \Res_\Gamma^\Delta \Ind_{\Gamma_i}^\Gamma V&\to&
\displaystyle{\bigoplus_{i\in I}} \Res_\Gamma^\Delta \Ind_{\Gamma_i}^\Gamma \Res_\Gamma^{\Gamma_i}V
&\to & \Res_\Gamma^\Delta V \\
\downarrow &&&& \downarrow \\
\displaystyle{\bigoplus_{j\in J}^{\ }} \Ind_{\Delta_j}^\Delta \Res_{\Gamma_i}^{\Delta_j} V &
\to& \displaystyle{\bigoplus_{j\in J}} \Ind_{\Delta_j}^\Delta V \ \to\ 
 \displaystyle{\bigoplus_{j\in J}} \Ind_{\Delta_j}^\Delta \Res_\Delta^{\Delta_j} V  &\to & V,\end{array}$$
 where the maps are defined as follows:
 \begin{itemize}
 \item The first downward arrow is defined by
 maps $\Res_\Gamma^\Delta \Ind_{\Gamma_i}^\Gamma V \to
 \Ind_{\Delta_j}^\Delta \Res_{\Gamma_i}^{\Delta_j} V$
 sending $f:\Gamma \to V$ to the map $\Delta \to V$ defined
 by $g\mapsto f(g_jw^{-1}gw)$.
 \item The final horizontal map in each row is defined by the evident trace map.
 \item The remaining maps are induced by the evident ones of the form $v \mapsto xv$
 where $x = w$, $w_j$, $y_i^{-1}$ or $z_j^{-1}$.
 \end{itemize}
 Choosing coset representatives $u_a \in \GL_2(F_{x_1})$ so that
 $U = \coprod_{a \in A} u_a U_1''$, decomposing
 $A = \coprod_{i\in I} A_i$ where
 $$A_i = \{\,a\in A\,|\, t u_a \in D^\times t_i U_1'' \A_{F,f}^\times\,\},$$
 and writing $tu_a = \delta_a t_i r_a$ with $\delta_a \in D^\times$, 
 $r_a \in U_1''\A_{F,f}^\times$ for each $a \in A_i$, we find that 
 $$\Gamma = \coprod_{a \in A_i} y_i^{-1}\Gamma_i y_i h_a$$
for each $i \in I$,  where $h_a = y_i^{-1} r_a u_a^{-1}$.
One can similarly choose coset representatives for each $\Delta_j$ in $\Delta$,
and the desired commutativity then follows from a direct calculation
using the resulting description of the trace maps as sums over $A$.
 \epf

\subsubsection{Badly dihedral representations}\

\begin{definition} \label{defn:bad} We say that $F'$ is a {\em $p$-bad quadratic extension}
of $F$ if $F'$ is a quadratic totally imaginary extension of $F$ of the form
$F(\delta)$ for some $\delta$ such that $\delta^p\in F^\times$ and
$\delta^p\CO_F = I^p\CO_F$ for some fractional ideal $I$ of $F$.
We say that an irreducible representation $\rho:G_F \to \GL_2(\Fpbar)$
is {\em badly dihedral} \footnote{There is a typo in the definition of {\em badly dihedral}  in the discussion before Lemma~4.11 of \cite{bdj}: $\delta^\ell \in K$ should be $\delta^\ell \in \CO_K$.
The definition here differs slightly from the one intended in \cite{bdj} in the case $p=2$
since we also wish to control the central character of the lift.}
if $\rho$ is induced from a character 
$G_{F'}\to \Fpbar^\times$ for some $p$-bad quadratic extension $F'$ of $F$.
\end{definition}

\begin{remark}  \label{rem:bad} Note that $F$ has a $p$-bad quadratic extension
if and only if $F$ contains the maximal real subfield of $F(\zeta_p)$.  If this is the case and $p$
is odd, then the only $p$-bad quadratic extension of $F$ is $F(\zeta_p)$, but if $p=2$,
then there are still only finitely many such extensions, as follows for example from the fact
that such an extension is necessarily unramified outside the primes dividing $2$ and $\infty$.
\end{remark}

Let $\rho:G_F\to\GL_2(\Fpbar)$ be the modular Galois representation from Theorem \ref{main}. Assume that $\rho$ arises from the definite quaternion algebra $D$ split at all finite places of $F$, in level $U$ and weight $V$ (here $V$ is a possibly reducible finite dimensional $\Fpbar$-linear representation of $\GL_2(\CO_F/(p))$). We keep the assumptions and notation from the previous section, so that in particular $\gm_\rho$ is the ideal of ${\bf T}^{S}_{\CO_E}$ attached to $\rho$, where $S$ is a finite set containing the primes of $F$ dividing $p$
and the primes at which $\rho$ is ramified.

Let $F_1,F_2,\ldots,F_r$ denote the $p$-bad quadratic extensions of $F$ (so $r \le 1$ 
unless $p=2$), and for each $i = 1,\ldots,r$, let $J_i$ denote the ideal of ${\bf T}^S_{\CO_E}$
generated by the elements $T_x$ for those finite places $x$ of $F$ such that $x\notin S$ and $x$ is inert in $F_i$.

\begin{lemma} \label{lem:just_inerts}  If $J_1J_2\cdots J_r \subset \gm_\rho$ then
$\rho$ is badly dihedral.
\end{lemma}
\begpf Since $\gm_\rho$ is prime, we may assume that $J_i \subset \gm_\rho$
for some $i$.  We thus have that $\tr(\rho(\Frob_v)) = 0$ for all $v\not\in S$ inert in $F_i$.
By the Chebotarev Density Theorem, it follows that $\tr(\rho(g)) =0$ for all
$g \in G_F \setminus G_{F_i}$.   Let $L$ denote the projective splitting field
of $\rho$, i.e., the fixed field of the kernel of the composite of $\rho$ with
the projection to $\PGL_2(\Fpbar)$.  

We claim that $F_i \subset L$.
Indeed if not, then we may choose $g \in G_L \setminus G_{F_i}$
and observe that for any $h \in G_F$, we have that either $h \not\in G_{F_i}$
so that $\tr(\rho(h)) = 0$, or $gh \not\in G_{F_i}$ in which case 
$\tr(\rho(g)\rho(h)) = \tr(\rho(gh)) = 0$ also implies that $\tr(\rho(h)) = 0$
since $\rho(g)$ is a scalar.  If $p > 2$, then taking $h$ to be the identity
immediately gives a contradiction; if $p=2$, then we see that every
element of $\gal(L/F)$ has order dividing $2$, which contradicts
the irreducibility of $\rho$.

It follows that $H = \gal(L/F_i)$ is subgroup of index $2$ in $G = \gal(L/F)$, and
that every element of $G\setminus H$ has order $2$.  Moreover $G$ is isomorphic
to a finite subgroup of $\PGL_2(\Fpbar)$ which is not contained in a Borel subgroup.
By Dickson's classification of such subgroups, we see the only possibility
is that $G$ is isomorphic to a dihedral group and $H$ is a cyclic subgroup
of index $2$.   Therefore the projective image of $\rho(G_{F_i})$ is cyclic,
from which it follows that $\rho|_{G_{F_i}}$ is reducible, and hence that $\rho$
is induced from a character of $G_{F_i}$.
\epf

\begin{lemma}\label{lem:support}  There is a finite set of places $S'$ such that
if $x_\nu \not\in S'$ and $x_\nu$ is inert in $F_\nu$ for $\nu=1,\ldots,r$, then
$T_{x_1}\cdots T_{x_r}$ annihilates $H^1(X_U,V)$.
\end{lemma}
\begpf  Write $D_f^\times = \coprod_i D^\times t_i U\A_{F,f}^\times$ and
choose a representative $t_i^{-1} \gamma t_i$ with $\gamma \in D^\times$
for each conjugacy class of elements of order $p$ in each of the groups $\Gamma_i=F^\times\backslash(U\A_{F,f}^\times\cap t_i^{-1}D^\times t_i)$.  Then $F[\gamma]$
is $p$-bad, so $F[\gamma] = F_\nu$ for some $\nu$.  Let $S_\gamma$ be the finite
set of places $x$ of $F$ inert in $F_\nu$ such that $\CO_{F_x}[\gamma] \neq \CO_{F_\nu,x}$,
and let $S'$ contain the union of the $S_\gamma$ for all $\gamma$ as above.

Now let $x_1,\ldots,x_r$ be as in the statement of the lemma and
let $T = T_{x_1\cdots x_r}$, which by Lemma~\ref{lem:product}
coincides with $T_{x_1}\cdots T_{x_r}$.  Let $U'$, $t_{ij}$ and $\Gamma_{ij}'$ be as
in the definition of the Hecke operator $T$ on $\oplus_i H^1(\Gamma_i,V)$ (cf. \ref{HeckeH1}).
We claim that $\Gamma_{ij}'$ has order prime to $p$.
Indeed if $t_{ij}^{-1}\gamma' t_{ij}$ is a representative of an element of
order $p$ in $\Gamma_{ij}'$, then its image in $\Gamma_i$ is of the form
$t_i^{-1}\gamma t_i$ for some $\gamma$ as above, so $\gamma$ is
conjugate in $\GL_2(F_x)$ to an element of $U'_xF_x^\times$,
where  $x = x_\nu$ for $\nu$ chosen so that $F[\gamma] = F_\nu$, and 
$$U_x' = U_0(x) = \left\{\left.\left(\begin{smallmatrix} a &b \\ c & d \end{smallmatrix}\right)
 \in \GL_2(\CO_{F,x})\,\right|\, c \equiv 0 \bmod x\right\}.$$
 Since $x\not\in S_\gamma$, we see that $\det\gamma \in \CO_{F,x}^\times$,
 so in fact $\gamma$ is conjugate
to an element of $U_0(x)$, hence its characteristic polynomial
is reducible mod $x$.  On the other hand since $x\not\in S_\gamma$,
we see also that $\gamma$ generates the ring of integers of an unramified
quadratic extension of $F_x$, so the characteristic polynomial of $\gamma$
is irreducible mod $x$, giving a contradiction.

Since all the $\Gamma'_{ij}$ have order prime to $p$, it follows that
$\oplus_{i,j} H^1(\Gamma'_{ij},V) = 0$, and therefore $T=0$.
\epf

Lemma~\ref{lem:exactness} follows easily from Lemmas~\ref{lem:just_inerts}
and~\ref{lem:support}.  Indeed it suffices to prove that if $\gm_\rho$ is in the
support of $H^1(X_U,V)$, then $\rho$ is badly dihedral.  Note that we may
enlarge $S$ since if $\gm_\rho$ is in the support of $H^1(X_U,V)$, then so
is $\gm'_\rho = \gm_\rho \cap {\bf T}^{S'}$ for any finite $S' \supset S$.
In particular choosing $S'$ as in Lemma~\ref{lem:support} we see that
the ideal $J_1\cdots J_r$ of Lemma~\ref{lem:just_inerts} is contained in the
annihilator of $H^1(X_U,V)$, so if $H^1(X_U,V)_{\gm_\rho} \neq 0$,
then $\rho$ is badly dihedral.

\subsection{The refinement of the argument in the case $[F:\Q]$ odd}\label{Second_odd}

Assume now that $[F:\Q]>1$ is odd, and let $\rho$ be as in the statement of Theorem \ref{main}.

Fix an infinite place $\tau_0$ of $F$ and let $D$ denote the quaternion algebra over $F$ whose ramification set equals $\Sigma-\{\tau_0\}$. Fix moreover isomorphisms $D_f\cong M_2(\A_{F,f})$ and $D\otimes_{F,\tau_0}\R\cong M_2(\R)$. Let $U=\prod_x U_x$ be an open compact subgroup of $D_f^\times\cong \GL_2(\A_f)$. Set $\h^{\pm}=\mathbf{P}^1(\C)-\mathbf{P}^1(\R)$. We denote by $X_U$ the Shimura curve over $F$ which is the canonical model for the complex analytic space $X_U(\C)=D^\times\backslash (D_f^\times\times \h^{\pm})/U\A_{F,f}^\times$  as in \cite{ca1} (we see $F\subset\C$ via $\tau_0$). Notice that we quotient $D_f^\times$ also by the action of  $\A_{F,f}^\times$ in order to keep track of central characters in what follows.

We have the decomposition into connected components:
$$X_U(\C)=\bigsqcup_{i\in I} \Gamma_i\backslash\h,$$
where $I$ is a finite set, and $\Gamma_i=F^\times\backslash \tilde\Gamma_i$ acts properly discontinuously on $\h$. Here $\tilde\Gamma_i=U\A_{F,f}^\times\cap t_i^{-1} D_+^\times t_i$  where
$D_f^\times = \coprod_{i\in I} D_+^\times t_i U\A_{F,f}^\times$ and $D_+^\times$ is the set of elements of $D$ of totally positive reduced norm.
The groups $\Gamma_i$ are torsion free and act freely and properly on $\h$ if $U$ is small enough, but not in general. Each component $\Gamma_i\backslash \h$ has a canonical model defined over a finite abelian extension of $F$ by \cite[1.2]{ca1}.

Let $S$ denote a finite set of finite places of $F$ containing the places above $p$ and the places at which $\rho$ is ramified. Let $\vec{k}\in\Z_{\geq 2}^\Sigma$ and $w\in\Z$ such that $k_\tau\equiv w \mod 2$ for all $\tau\in\Sigma$ and let $\psi$ be a finite order Hecke character of $\A_F^\times$ totally of parity $w$. The system of Hecke eigenvalues for the action of $\T_\C^S$ on the space of holomorphic Hilbert modular forms of level $U$, weight $(\vec{k},w)$, and central character $\psi^{-1}|\cdot|^{-w}$ coincide with the systems of Hecke eigenvalues arising from the \'etale cohomology $H^1(X_{U,\overline{F}},V_{\vec{k},w,\psi}\otimes_{\CO_E}E)$, where $V_{\vec{k},w,\psi}$ is the $\CO_E$-sheaf associated to the homonymous representation of $U_p\A_{F,f}^\times$ on a finite free $\CO_E$-module (see \ref{Sub:DQA} and \cite{ca}). Here $E$ is a large enough finite extension of ${\Q}_p$, and we see $E\subset\overline{\Q}_p\cong\C$. These Hecke eigensystems also coincide with the Hecke eigensystems arising from  $\oplus_{i\in I} H^1(\Gamma_i,V_{\vec{k},w,\psi}\otimes_{\CO_E}E)$. (The action of the Hecke operators on the cohomology of the $\Gamma_i$'s is  defined as in section \ref{HeckeH1}.) We denote by $\gm_\rho$ the prime ideal of $\T^S_{\CO_E}$ attached to $\rho$.

To prove the second half of Theorem \ref{main} we may use the same strategy adopted in the case of a definite quaternion algebra, modulo guaranteeing that an analogue of Lemma \ref{lem:exactness} holds.  We now consider the natural action of $\T^S_{\CO_E}$ on spaces $(\oplus_{i\in I}H^j(\Gamma_i,V))_{\gm_\rho}$ where $V$ is a finite dimensional continuous $\Fpbar$-linear representation of $U_p\A_{F,f}^\times/F^\times$; we wish to prove that if $\rho$ is not badly dihedral,
then the functor $V \mapsto (\oplus_{i\in I}H^1(\Gamma_i,V))_{\gm_\rho}$ is exact, so it is enough to prove
that  $(\oplus_{i\in I}H^j(\Gamma_i,V))_{\gm_\rho} = 0$ for $j=0,2$.

First note that by the strong approximation theorem, the reduced norm $\det$ induces
a bijection $D_+^\times\backslash D_f^\times /U\A_{F,f}^\times \to I$
where $I = \A_F^\times / F_{\infty,>0}^\times F^\times \det(U) (\A_{F,f}^\times)^2$.
It follows that in the definition of
the Hecke operator $T_x$ for places $x \not\in S$, we may use the same index set
$I$ and representatives $t_i$ when $U$ is replaced by $U' = U \cap \Pi_x^{-1}U\Pi_x$ and
$U'' = \Pi_x U' \Pi_x^{-1}$.  We may thus write $T_x$ as a direct sum of composite maps
$$H^j(\Gamma_i,V) \to H^j(\Gamma'_i,V) \to H^j(\Gamma''_{i'},V) \to H^j(\Gamma_{i'},V)$$
where $i \mapsto i'$ is induced by multiplication by $(\varpi_x)_x$ on $I$.  Let $F_I$
denote the abelian extension of $F$ corresponding to $I$ by class field theory; we see that
if $x$ splits completely in $F_I$, then $T_x$ acts componentwise on
$\oplus_{i\in I} H^j(\Gamma_i,V)$.
Moreover if $j = 0$, then for such $x$ this action is simply multiplication by the index
$[\Gamma_i:\Gamma_i'] = [\Gamma_i:\Gamma_i''] = \N(x) + 1$ on each component.

Suppose now that $\gm_\rho$ is in the support of $\oplus_{i\in I} H^0(\Gamma_i,V)$,
and let $F'$ denote  the composite of $F_I(\zeta_p)$ with the splitting field of $\det(\rho)$.
If $x$ splits completely in $F'$, then $T_x - \N(x) - 1 \in \gm_\rho$ for all such $x$, which
implies that $\tr(\rho(\Frob_x)) = 2$.   The Brauer-Nesbitt and Chebotarev Density Theorems
then imply that $\rho|_{G_{F'}}$ has trivial semi-simplification; since $F'$ is an abelian 
extension of $F$, this contradicts the irreducibility of $\rho$.

To treat the case of $j=2$, we use Farrell cohomology groups $\widehat{H}^j(\Gamma,V)$
(defined in \cite{fa}) for finite index subgroups $\Gamma$ of the groups $\Gamma_i$.
Note that if $F = \Q$, then such $\Gamma$
have virtual cohomological dimension one, so that $H^2(\Gamma,V) = \widehat{H}^2(\Gamma,V)$.
If $F\neq \Q$, then $\Gamma$ is a virtual duality group of dimension $2$ with dualizing
module isomorphic to $\Z$ (with trivial $\Gamma$ action)\footnote{Note that there is a canonical choice of orientation
$H^2(\Gamma,\Z\Gamma)\cong \Z$ provided by the complex analytic structure on $X_U(\C)$.}, so that \cite[Thm.~2]{fa}
yields an exact sequence:
$$H_0(\Gamma,V)  \to H^2(\Gamma,V)\to \widehat H^2(\Gamma,V)\to 0.$$
We will assume $F\neq \Q$ since the case $F = \Q$ is easier and can be
treated by minor modifications to the arguments below.

For $\Gamma'$ a finite index subgroup of $\Gamma$, we have restriction
maps $H_j(\Gamma,V) \to H_j(\Gamma',V)$ and $\widehat{H}^j(\Gamma, V) \to \widehat{H}^j(\Gamma',V)$, 
as well as corestriction maps $H_j(\Gamma',V) \to H_j(\Gamma,V)$ and
$\widehat{H}^j(\Gamma', V) \to \widehat{H}^j(\Gamma,V)$, allowing us to define Hecke operators
$T_x$  on $\oplus_{i\in I} H_j(\Gamma_i,V)$ and $\oplus_{i\in I} \widehat{H}^j(\Gamma_i,V)$ 
for $x\not\in S$ exactly as on  $\oplus_{i\in I} H^j(\Gamma_i,V)$.   By the following lemma
(and the fact that the isomorphisms $\Gamma_i' \cong \Gamma_i''$ are orientation-preserving),
the homomorphisms 
\begin{equation}\label{H2seq}
\bigoplus_{i\in I} H_0(\Gamma_{i},V)
\to \bigoplus_{i\in I} H^2(\Gamma_{i},V)\to \bigoplus_{i\in I} \widehat H^2(\Gamma_{i},V)
\end{equation}
are compatible with the operators $T_x$.  

\begin{lemma}
Let $\Gamma$ be a virtual duality group of dimension $n$ with dualizing module $D$.
Let $M$ be a left $\Z\Gamma$-module
and $\Gamma'$ a finite index subgroup of $\Gamma$.  Then the diagram:
\small
$$\begin{array}{ccccccccc}
\cdots\rightarrow\!\!\!\!\!&H_{n-j}(\Gamma,D\otimes_{\Z}M)&\!\!\!\!\!\rightarrow\!\!\!\!\!&H^j(\Gamma,M)&\!\!\!\!\!\rightarrow\!\!\!\!\!&
\widehat{H}^j(\Gamma,M)&\!\!\!\!\!\rightarrow\!\!\!\!\!&H_{n-j-1}(\Gamma,D\otimes_{\Z}M)&\!\!\!\!\!\rightarrow\cdots\\
&\downarrow&&\downarrow&&\downarrow&&\downarrow&\\
\cdots\rightarrow\!\!\!\!\!&H_{n-j}(\Gamma',D\otimes_{\Z}M)&\!\!\!\!\!\rightarrow\!\!\!\!\!&H^j(\Gamma',M)&\!\!\!\!\!\rightarrow\!\!\!\!\!&
\widehat{H}^j(\Gamma',M)&\!\!\!\!\!\rightarrow\!\!\!\!\!&H_{n-j-1}(\Gamma',D\otimes_{\Z}M)&\!\!\!\!\!\rightarrow\cdots\\
\end{array}$$
\normalsize
commutes, where the rows are the exact sequences given by \cite[Thm.~2]{fa}
and the vertical arrows are the natural restriction maps.   Similarly the diagram
commutes with the downward arrows replaced by the upward corestriction maps.
\end{lemma}
\begpf  Let $(P_\bullet,d_\bullet)$ be a projective resolution of finite type of the trivial (left\footnote{Some of the modules we consider will be naturally right $\Z\Gamma$-modules; they can be regarded as left $\Z\Gamma$-modules via the involution $\gamma\mapsto \gamma^{-1}$ of $\Gamma$; and vice versa. Some of the chain complexes we consider will be sometimes regarded as cochain complexes, after relabelling; and vice versa.})
$\Gamma$-module $\Z$, which we view also as a projective resolution of $\Z$ as a 
$\Gamma'$-module.  We define $P_\bullet^\ast:=\mathrm{Hom}_\Gamma(P_\bullet,\Z\Gamma)$ and ${P'_\bullet}^\ast:=\mathrm{Hom}_{\Gamma'}(P_\bullet,\Z\Gamma')$ and denote by $(P_\bullet^\ast,d_\bullet^\ast)$ and $({P'_\bullet}^\ast,{d'_\bullet}^\ast)$ the corresponding \emph{cochain} complexes of right $\Z \Gamma$- and $\Z\Gamma'$-modules respectively. There is a natural map of cochain complexes of right $\Z\Gamma'$-modules $\rho_\bullet:P_\bullet^\ast\to {P_\bullet'}^\ast$ induced by the map $ \Z\Gamma\to\Z\Gamma'$ given by $\sum_{\gamma\in\Gamma}n_\gamma \gamma \mapsto \sum_{\gamma\in\Gamma'} n_\gamma \gamma$. 
Note that $\rho_\bullet$ is an isomorphism, with inverse $\sigma_\bullet$ defined
by $(\sigma_j(f))(x) = \sum_{\gamma\in \mathscr{R}} \gamma^{-1}f(\gamma x)$
for $f\in {P'_j}^\ast$ and $x\in P_j$, where $\Gamma=\sqcup_{\gamma\in\mathscr{R}}\Gamma'\gamma$.

Recall that the dualizing module $D$ is defined as $H^n(\Gamma,\Z\Gamma)$, which
we view as a right $\Z\Gamma$-module, and let $(Q_\bullet,e_\bullet)$ be a projective resolution of $D$ as a right $\Z\Gamma$-module. Note that $D=H^n(P_\bullet^\ast)=\ker d_n^\ast/\mathrm{Im\,} d_{n-1}^\ast$, and then the natural inclusion $D\hookrightarrow \mathrm{coker\,} d_{n-1}^\ast$ can be extended to a map of \emph{chain} complexes $f_\bullet:Q_\bullet \to P_{n-\bullet}^\ast$.
Moreover if we let $D' = H^n(\Gamma',\Z\Gamma')$ denote the dualizing module of $\Gamma'$,
then $\rho_\bullet$ induces the canonical isomorphism $D \cong D'$ of $\Z\Gamma'$-modules,
so that we may also view $(Q_\bullet,e_\bullet)$ as a projective resolution of $D'$, and extend
the natural inclusion $D'\hookrightarrow \mathrm{coker\,}(d_{n}')^\ast$ to a map of chain complexes $f'_\bullet:Q_\bullet \to (P_{n-\bullet}')^\ast$ where $f_\bullet' = \rho_\bullet \circ f_\bullet$. 

We now let $X_\bullet$ denote the mapping cone of the chain map $f_\bullet$, so that $X_\bullet=Q_\bullet\oplus P_{n-\bullet-1}^\ast$, and similarly let $X_\bullet'$ be the mapping cone of $f'_\bullet$.
Then $\mathrm{id} \oplus \rho_\bullet$ defines a chain map, giving a commutative diagram of morphisms
\emph{cochain} complexes of right $\Z\Gamma'$-modules: 
\begin{equation}\label{phi} 
\begin{array}{ccccccccc}
0 & \to & P_\bullet^\ast & \to & X_{n-\bullet-1} & \to & Q_{n-\bullet -1} & \to & 0  \\
  &  &  \downarrow & & \downarrow & & \downarrow & &\\
0 & \to & {P_\bullet'}^\ast & \to & X'_{n-\bullet-1} & \to & Q_{n-\bullet -1} & \to & 0 
\end{array} 
\end{equation}
in which the rows are exact and the vertical maps are isomorphisms.

We now apply the functor $(\,\cdot\,)\otimes_\Gamma M$ to the first line of \eqref{phi}, and the functor $(\,\cdot\,)\otimes_{\Gamma'} M$ to the second.  For a right $\Gamma$-module $A$ and
 a left $\Gamma$-module $B$, we define the trace map 
 $\tr: A \otimes_\Gamma B \to A\otimes_{\Gamma'} B$
 by $\tr(a\otimes b) =  \sum_{\gamma\in\mathscr{R}} a\gamma^{-1}\otimes\gamma b$.
We thus obtain a commutative diagram of complexes with exact rows:
\begin{equation}\label{phi_M} 
\begin{array}{ccccccccc}
0 & \to & P_\bullet^\ast\otimes_\Gamma M & \to & X_{n-\bullet-1}\otimes_\Gamma M & \to & Q_{n-\bullet -1}\otimes_\Gamma M & \to & 0  \\
  &  &  \downarrow & & \downarrow & & \downarrow & &\\
0 & \to & {P_\bullet'}^\ast \otimes_{\Gamma'} M& \to & X'_{n-\bullet-1} \otimes_{\Gamma'} M& \to & Q_{n-\bullet -1} \otimes_{\Gamma'} M& \to & 0 
\end{array} 
\end{equation}
where the left vertical arrow is given by
$P_\bullet^\ast\otimes_\Gamma M \overset{\mathrm{tr}}\longrightarrow
P_\bullet^\ast\otimes_{\Gamma'} M \overset{\rho_\bullet\otimes\mathrm{id}_M}\longrightarrow
{P_\bullet'}^\ast\otimes_{\Gamma'}M$,
the right vertical arrow is the trace map, and the middle vertical arrow
is their direct sum.   

Taking cohomology in \eqref{phi_M} then yields the desired
commutative diagram.  Indeed the long exact sequences of \cite[Thm.~2]{fa} are
precisely those associated to the rows of \eqref{phi_M}, and it is straightforward
to check that the vertical maps induce the corresponding restriction maps on
homology and cohomology (for Farrell cohomology, this follows from the
characterization of $\res$ following \cite[Rem.~2]{fa}).

The proof of compatibility with corestriction is similar, so we omit the details.
One just uses $\sigma_\bullet$ instead of $\rho_\bullet$ to obtain a diagram
as in \eqref{phi}, but with upward arrows, and use the canonical projection
$A\otimes_{\Gamma'}B \to A\otimes_\Gamma B$ to obtain the analogue
of \eqref{phi_M}, again with upward arrows.
\epf

We can now use \eqref{H2seq} to prove that $\gm_\rho$ is not in the
support of $\oplus_{i\in I} H^2(\Gamma_{i},V)$. Indeed
the same argument as for $H^0$ shows that if $x$ splits completely in $F_I$,
then $T_x = \N(x) + 1$ on $\oplus_{i\in I} H_0(\Gamma_{i},V)$, so that the
irreducibility of $\rho$ implies that $\gm_\rho$ is not in the support of the image
of $\oplus_{i\in I} H_0(\Gamma_{i},V)$.  Note also that the surjectivity of  
$\oplus_{i\in I} H^2(\Gamma_{i},V)\to \oplus_{i\in I} \widehat H^2(\Gamma_{i},V)$
implies that the operators $T_x$ commute, hence $\T^S_{\CO_E}$ acts,
on $\oplus_{i\in I} \widehat H^2(\Gamma_{i},V)$.  Thus it suffices to prove
that if $\rho$ is not badly dihedral, then $\gm_\rho$ is not in the support
of $\oplus_{i\in I} \widehat H^2(\Gamma_{i},V)$.

Let $S'$ be a finite set of finite places of $F$ constructed as in the proof of Lemma \ref{lem:support}.
(Now $D$ is an indefinite quaternion algebra, so the groups $\Gamma_i$ are infinite, but each
still has only finitely many conjugacy classes of elements of order $p$.)
For each $\nu=1,\dots,r$ let $x_\nu\not\in S'\cup S$ be a finite place of $F$ inert in $F_\nu$;
let $T = T_{x_1\cdots x_r} = T_{x_1}\cdots T_{x_r}$, and let $\Gamma_i'$ be as in
the definition of the Hecke operator $T$.  (Note that we can use the same index set $I$
and representatives $t_i$.)   The same proof as in Lemma \ref{lem:support} shows that
$\Gamma_i'$ does not contain any element of order $p$.   Therefore $\Gamma_i'$
has a torsion-free subgroup of finite index prime to $p$, so $\widehat{H}^2(\Gamma_i',V) = 0$.
It follows that the operator $T$ annihilates $\oplus_{i\in I} \widehat{H}^2(\Gamma_i,V)$,
since it factors through $\oplus_{i \in I} \widehat{H}^2(\Gamma_i',V)$.  To prove that 
$(\oplus_{i\in I}\widehat{H}^2(\Gamma_i,V))_{\gm_\rho}=0$, we may enlarge $S$ so that $S \supset S'$. For each $\nu=1,\dots,r$ we denote by $J_\nu$ the ideal of $\T_{\CO_E}^S$ generated by the Hecke operators $T_x$ for those finite places $x$ of $F$ such that $x\not\in S$ and $x$ is inert in $F_\nu$. Observe that the ideal $J_1J_2\cdots J_r$ annihilates $(\oplus_{i\in I}\widehat{H}^2(\Gamma_i,V))_{\gm_\rho}$. By Lemma \ref{lem:just_inerts} (which holds, \emph{mut. mut.}, also in the current setting) we deduce that $(\oplus_{i\in I}\widehat{H}^2(\Gamma_i,V))_{\gm_\rho}$ vanishes, since $\rho$ is not badly dihedral.
This completes the proof that $(\oplus_{i\in I}H^2(\Gamma_i,V))_{\gm_\rho} = 0$, and hence the
functor $V \mapsto (\oplus_{i\in I} H^1(\Gamma_i,V))_{\gm_\rho}$ is exact.

\section{A variant in a special case}

We now give a variant of the main result in a special case, with a view
to producing forms satisfying  the hypotheses of Assumption~8.15 in Section~8.3 of \cite{dp}.

We must first introduce some notation.  Recall that we have fixed an embedding
$\ol{\Q} \to \C$ and an isomorphism $\ol{\Q}_p \cong \C$.   Note that these choices induce
a bijection between the set $\Sigma$ of embeddings $\tau: F\to \R$ and the set of pairs
$(v,\vartheta)$ where $v|p$ and $\vartheta$ is an embedding $F_v \to \ol{\Q}_p$.
For each $v|p$ we let $\Sigma_v$ denote the set of embeddings $\vartheta: F_v \to \ol{\Q}_p$,
which we identify with a subset of $\Sigma$ via this bijection. 

Let $\pi$ be a cuspidal automorphic representation of $\GL_2(\A_F)$ which is
holomorphic of weight $(\vec{k},w)$, where as usual $\vec{k} \in \Z_{\ge 2}^\Sigma$
and $w \in \Z$ is such that $w \equiv k_\tau \bmod 2$ for all $\tau \in \Sigma$.
Suppose further that for all $v|p$, the local factor $\pi_v$ is either unramified
principal series or an unramified twist of the Steinberg representation.
Let $a_v(\pi)$ denote the eigenvalue of the Hecke operator
$T_v = U_v \Pi_v U_v$ on the one-dimensional vector space $\pi_v^{U_v}$,
where $\Pi_v=\left(\begin{smallmatrix} \varpi_v&0\\ 0&1 \end{smallmatrix}\right)\in \GL_2({F_v})$
and $U_v = \GL_2(\CO_{F,v})$ or $U_0(v)$ according to whether $\pi_v$ is
unramified.
Since $a_v(\pi)$ is algebraic, we may view it as an element of $\ol{\Q}$ via
our choices of embeddings $\ol{\Q} \to \C$ and $\ol{\Q} \to \ol{\Q}_p$.
We say that $\pi$ is {\em ordinary} at $v$ (with respect to our choices of
embeddings) if
$$|a_v(\pi)|_p = p^{\sum_{\tau \in \Sigma_v} (k_\tau - 2- w)/2e_v}.$$
\begin{remark} In general we have the expression on the right as an upper bound
on $|a_v(\pi)|_p$; this follows for example from \cite[Thm.~4.11]{hi}, but will also
be clear from the proof of Theorem~\ref{thm:variant1} below.   Moreover if 
equality holds then Theorem~1 of \cite{sk} implies that
the local Galois representation $\rho_\pi|_{G_{F_v}}$
is reducible.  Note that if $\pi_v$ is an unramified twist of the Steinberg
representation, then $\pi$ is ordinary at $v$ if and only if $k_\tau = 2$
for all $\tau \in \Sigma_v$.
\end{remark}

\begin{theorem}  \label{thm:variant1}
Suppose that $\rho:G_F \to \GL_2(\Fpbar)$ is such that 
$\rho \cong \ol{\rho}_\pi$ for some cuspidal, holomorphic, automorphic
representation $\pi$ of $\GL_2(\A_F)$ of weight $(\vec{k},w) = (2,\dots,2)$
such that for each $v|p$, $\pi_v$ is either unramified principal series or
an unramified twist of the Steinberg representation.  For any finite set of
primes $T$ of $F$, there exist a  cuspidal automorphic representation $\pi'$
of $\GL_2(\A_F)$ and a character $\xi:G_F \to \ol{\F}_p^\times$ of order at most $2$
such that
\begin{itemize}
\item if $\tau \in \Sigma$, then $\pi_\tau' \cong D_{k'_\tau,w'}$ with $k'_\tau \in \{2,w'\}$
  where $w' = p+1$ if $p$ is odd and $w' = 4$ if $p=2$;
\item if $v|p$, then $\pi'_v$ is unramified principal series, and is ordinary if $k'_\tau = w'$
   for some $\tau \in \Sigma_v$;
\item the prime-to-$p$ part of the conductor of $\xi$ divides a prime
   $y \not\in T$ which splits completely in $F$;
\item $\ol{\rho}_{\pi'} \cong \xi \otimes \rho$.
\end{itemize}
Suppose further that $\pi$ has prime-to-$p$ conductor dividing $\gn \subset \CO_F$,
and that $\psi$ is a totally even finite order Hecke character of $\A_F^\times$ of conductor
dividing $\gn$ satisfying $\det\rho = \overline{\psi\epsilon}$.  Then if $\rho$ is not badly
dihedral, we can choose $\pi'$ as above with conductor dividing $\gn y^2$, central
character $\psi^{-1}|\ |^{-w'}$ and $\xi_y \otimes \pi'_y$ unramified principal series.
\end{theorem}

\begin{remark}  We will see from the proof that  the conclusion can be made more
precise as follows:  For each $v|p$ such that $\pi_v$ is ramified, we can ensure
that $\pi'_v$ is ordinary and the set of $\tau \in \Sigma_v$ such that $k'_\tau = w'$
maps bijectively to $\ol{\Sigma}_v$ under the natural projection.
\end{remark}

\begpf   Let $R$ denote the set of primes $v|p$ such that $\pi_v$ is ramified,
and as usual let $S$ be a sufficiently large finite set of primes containing all those
dividing $p$ and all those at which $\pi$ is ramified.
We suppose that $E$ is a sufficiently large finite extension of $\Q_p$ in $\ol{\Q}_p$
that contains the eigenvalue $a_x(\pi)$ of $T_x$ on $\pi_x^{\GL_2(\CO_{F,x})}$ for all
$x \not\in S$ and (necessarily also) the eigenvalue $a_v(\pi)$ of $T_v$ on $\pi_v^{U_0(v)}$
for all $v\in R$.

Let $S' = S\setminus R$ and let $\T = \T_{\CO_E}^{S'}$ denote the $\CO_E$-algebra
generated by the operators $T_x$ for $x \not\in S'$.    Let $\gm$ denote the kernel
of the homomorphism $\T \to \Fpbar$ defined by sending $T_x$ to the reduction of
of $a_x(\pi\otimes|\det|)) = a_x(\pi)\N(x)^{-1}$ for $x \not\in S'$.
(For convenience in keeping track of ordinariness, we have replaced $\pi$ by its
twist by $|\det|$  to ensure that $T_v\not\in \gm$ for $v \in R$.)

Let $U = U_1(\gn) \cap U_0(\prod_{v\in R}v)$ where $\gn$ is the prime-to-$p$ conductor of $\pi$.
If $[F:\Q]$ is even, then we let $D$ be the definite quaternion algebra over $F$ ramified at precisely the infinite places of $F$.
By the Jacquet--Langlands correspondence, $\gm$ is in the support of $S_{V'}(U)$
where $V'$ is the representation of  $U\A_{F,f}^\times$ on $\Fpbar$ on which $U$ acts
trivially and $\A_{F,f}^\times$ acts via $\ol{\psi}$.  If $[F:\Q]$ is odd, then we let $D$
be a quaternion algebra over $F$ ramified at precisely all but one infinite place.
In this case $\gm$ is in the support of $\oplus_{i\in I} H^1(\Gamma_i,V')$ where
as before, $D_f^\times = \coprod_{i\in I} D_+^\times t_i U\A_{F,f}^\times$ and
$\Gamma_i=F^\times\backslash (U\A_{F,f}^\times\cap t_i^{-1} D_+^\times t_i)$.  Since the
argument is the same in the case of either parity, we will ease notation by
writing $S(U,V')$ for $\oplus_{i\in I} H^j(\Gamma_i,V')$ where $j = 0$ or $1$
according to the parity of $[F:\Q]$ (so $S(U,V') = S_{V'}(U)$ if $j=0$).

We will now show that $\gm$ is in the support of $S(U_1(\gn),V)$ where
$V = V' \otimes (\otimes_{v\in R}S_{v,(p-1,\ldots,p-1)})$.  We proceed by induction
on $|R'|$ to show that $\gm$ is in the support of $S(U_{R'},V_{R'})$ for
$R' \subset R$, where $U_{R'} = U\prod_{v\in R'}\GL_2(\CO_{F,v})$
and $V_{R'} = V' \otimes (\otimes_{v\in R'}S_{v,(p-1,\ldots,p-1)})$.
Note that $U_R = U_1(\gn)$, $U_\emptyset = U$, and
we already know that $\gm$ is in the support of $S(U_\emptyset, V_\emptyset)$.

Suppose now that $v \in R \setminus R'$.
The canonical isomorphisms 
 $$S(U_{R'},V_{R'})\cong  S(U_{R'\cup\{v\}} ,\Ind_{U_{R'}}^{U_{R'\cup\{v\}}}(V_{R'}))$$
  and $\Ind_B^{\GL_2(k_v)}\ol{\F}_p \cong \ol{\F}_p \oplus S_{(p-1,\ldots,p-1)}$
give rise to an exact sequence
\begin{equation}\label{eqn:Gamma0} 0 \to S(U_{R'\cup\{v\}},V_{R'})
      \stackrel{\alpha_v}{\to} S(U_{R'},V_{R'}) 
      \stackrel{\beta_v}{\to}   S(U_{R'\cup\{v\}},V_{R'\cup\{v\}}) \to 0\end{equation}
such that $\alpha_v$ and $\beta_v$ are
compatible with the operators $T_x$ for $x \not\in S' \cup \{v\}$,
but not necessarily with the operator $T_v$.  Note however that
the matrix $w_v = \left(\begin{smallmatrix} 0&1\\ \varpi_v&0 \end{smallmatrix}\right)
 \in \GL_2(F_v)$ normalizes $U_{R'}$ so that $W_v = (w_v)_*$ defines an
automorphism of $S(U_{R'},V_{R'})$ which is compatible with
$T_x$ for $x\not\in S'\cup\{v\}$ and satisfies $W_v^2 = \ol{\psi}(\varpi_v)^{-1}$.
Moreover unravelling the definitions of the operator $T_v$ one finds that
$T_v W_v \alpha_v = 0$ and $\beta_v W_v T_v = T_v \beta W_v$.
Therefore $\beta_v W_v$ is $\T$-linear, and $\gm$ is not in the support
of its kernel since $T_v \not\in \gm$.   It follows that if $\gm$ is in the support of
$S(U_{R'},V_{R'})$ then it is also in the support of 
$S(U_{R'\cup\{v\}},V_{R'\cup\{v\}})$.

To set the stage for lifting to characteristic zero, we distinguish between
the cases $p=2$ and $p> 2$.    

If $p>2$, then we
let $\ga$ denote the ideal of $\CO_F$ such that $p\CO_F = \ga \prod_{v\in R} v$.
Let $y$ be any prime ideal of $\CO_F$ such that $y\not\in S$ and
 $[y] = [\ga]^{-1}$ in the ray class group
of conductor $4\CO_F$ (if non-trivial; otherwise we can let $y = \CO_F$), and
choose $\alpha$ such that $y\ga = \alpha\CO_F$ and
$\alpha \equiv 1 \bmod 4\CO_F$.   Finally let $\xi$ be the character of $G_F$ with
splitting field $F(\sqrt{\alpha})$.  Thus $\xi$ is ramified precisely at $y$ and at
certain $v|p$, namely those such that either $v\in R$ and $e_v$ is even, 
or $v\not\in R$ and $e_v$ is odd.  Note also that $y$ can be chosen to
split completely in $F$.

If $p = 2$, then we let $\xi$ be the trivial character, but we must make another
modification instead of introducing a quadratic twist.   For each $v\in R$ we
have a $\GL_2(k_v)$-equivariant inclusion $S_{(1,\ldots,1)} \to S_{(2,\ldots,2)}$,
and these induce a $\T$-equivariant map
$S(U_1(\gn), V) \to
  S(U_1(\gn), V' \otimes (\otimes_{v\in R}S_{v,(2,\ldots,2)}))$.
One checks as usual that $\gm$ is not in the support of the kernel, so we can
replace each $S_{v,(p-1,\ldots,p-1)}$ by $S_{v,(2,\ldots,2)}$ in the definition
of $V$ for $p=2$.

Now let $\tilde{R} \subset \Sigma = \coprod_{v|p} \Sigma_v$ be a set of
embeddings $F \to \ol{\Q}$ such
that the $\tilde{R} \cap \Sigma_v = \emptyset$ if $v\not\in R$ and the natural
map $\tilde{R} \cap \Sigma_v \to \ol{\Sigma}_v$ is bijective if $v \in R$.
Define $\vec{k}' \in \Z_{\ge 2}^\Sigma$ by setting $k'_\tau = 2$ if $\tau \not\in \tilde{R}$
and $k'_\tau = w'$ if $\tau \in \tilde{R}$.  (Recall that $w' = p+1$ if $p>2$ and
$w' = 4$ if $p=2$.)   Consider the representation $V_{\vec{k'},w'-2,\psi}$
of $U_1(\gn)\A_{F,f}^\times$; recall that
this is the free $\CO_E$-module defined by 
$$V_{\vec{k'},w'-2,\psi} =  \left(\bigotimes_{\tau\in \tilde{R}} \Sym^{w'-2}\CO_E^2\right)
          \bigotimes \left(\bigotimes_{\tau\not\in\tilde{R}} \det^{(w'-2)/2}\right)$$
as a representation of  $\GL_2(\CO_{F,p})$, with $x \in \A_{F,f}^\times$ acting via 
$\N(x_p)^{w'-2}|x|^{w'-2}\psi(x)$.
We let $V_\xi = V_{\vec{k'},w'-2,\psi}$ if $p=2$ (or if $y = \CO_F$); otherwise
we let $V_\xi$ denote the twist of $V_{\vec{k'},w'-2,\psi}$ by the Teichmuller lift
of the character $\xi\circ\det$ of $\GL_2(\CO_{F,y})$. 

One finds that if $p > 2$, then the reduction of 
$\prod_{\tau \in \Sigma_v \setminus \tilde{R}}\det^{(p-1)/2}$ is $\det^{e_v(p^{f_v}-1)/2}$
if $v \in R$ and $\det^{(e_v-1)(p^{f_v}-1)/2}$ if $v\not\in R$, from which it follows that
 $\ol{V}_\xi$ is isomorphic to the twist of $V$ by the character $\xi \circ \det$
 of $U_1(\gn)\A_{F,f}^\times$.  The isomorphisms $V \to \ol{V}_\xi$ of
 $\Gamma_i = F^\times\backslash (U_1(\gn)\A_{F,f}^\times \cap t_i^{-1}D^\times t_i)$--modules defined by $v \mapsto \xi(\det(t_i)) v$
 induce an isomorphism
 $S(U_1(\gn),V) \cong S(U_1(\gn),\ol{V}_\xi)$ under which
 $\xi(\varpi_v)T_v$ corresponds to $T_v$ for $v$ not dividing $yp$ and to
 $T_v^0$ for $v\in R$, where $T_v^0$ is compatible with an operator on $S(U,V_\xi)$
 such that $T_v = T_v^0 \prod_{\tau\in \Sigma_v\setminus\tilde{R}} \tau(\varpi_v)^{(w' - 2)/2}$.
 (Note that $T_v^0$ may depend on the choice of uniformizer $\varpi_v$.)
 
 Now let $\T'$ denote the $\CO_E$-algebra generated by the operators $T_x$ for
 $x \not\in S\cup\{y\}$ and $T_v^0$ for $v \in R$, and let $\gm'$ denote the kernel
 of the homomorphism $\T' \to \Fpbar$ sending each $T_x$ to the reduction of 
 $\xi(\varpi_x)a_x(\pi)\N(x)^{-1}$ and each $T_v^0$ to the reduction of 
 $\xi(\varpi_v)a_v(\pi)\N(v)^{-1}$.   We then have that $\gm'$ is in the support 
of  $S(U_1(\gn),\ol{V}_\xi)$.  If $\rho$ is not badly dihedral, it follows
as in the proof of Theorem~\ref{main} that $\gm'$ is in the support of 
$S(U_1(\gn),V_\xi \otimes_{\CO_E} E)$, and hence that there is an
automorphic representation whose twist by $|\det|^{-1}$ is the required $\pi'$.
If $\rho$ is badly dihedral, then the proof goes through after replacing $U_1(\gn)$
by a smaller open compact subgroup $U$ so that the groups $\Gamma_i$
are torsion-free.
\epf
 
 \begin{remark}  \label{lassina} The necessity of the quadratic twist $\xi$ in the conclusion
follows from consideration of the local Galois representations $\ol{\rho}_{\pi'}|_{G_{F_v}}$
for $v|p$.  Furthermore one can construct explicit examples showing that $\xi$ may
need to be ramified outside $p$; we are grateful to L.~Demb\'el\'e for providing the
following one.  Over $F = \Q(\alpha)$ with $\alpha = \sqrt{10}$, there is 
a Hilbert modular form\footnote{Demb\'el\'e also offers the equation $y^2 = x^3 + (990144\alpha + 3127248)x -545501952\alpha - 1726178688$ for the associated elliptic curve.} of weight $(2,2)$,
level $(\alpha+2)$ and trivial character with leading Hecke eigenvalues (ordered by norm):
\scriptsize
$$\begin{array}{c|cccccccc}
v& (2, \alpha) & (3,\alpha+2) &  (3, \alpha + 4) & 
    (5,\alpha)&  (13,\alpha+6) & (13, \alpha+7) &
    (31, \alpha + 17)& (31, \alpha + 14) 
\\ \hline
a_v&-1&-1&3&1&-7&0&-3&4
\end{array}$$
\normalsize 
The corresponding automorphic representation $\pi$
then satisfies the hypotheses of the theorem for $p=3$, but one can show that
the character $\xi$ in the conclusion must be ramified at $(3,\alpha+4)$ but not
$(3,\alpha+2)$, from which it follows that $\xi$ must also be ramified at a prime $y$ not dividing $3$;
in fact one can let $y$ be any non-principal prime of $\CO_F$ not dividing $6$.
\end{remark}

 We now give another variant which under the same hypotheses produces
 lifts of parallel weight without the quadratic twist, at the expense of ordinariness.

\begin{theorem}  \label{thm:variant2}
Suppose that $\rho:G_F \to \GL_2(\Fpbar)$ is such that 
$\rho \cong \ol{\rho}_\pi$ for some cuspidal, holomorphic, automorphic
representation $\pi$ of $\GL_2(\A_F)$ of weight $(\vec{k},w) = (2,\dots,2)$
such that for each $v|p$, $\pi_v$ is either unramified principal series or
an unramified twist of the Steinberg representation.   Let $k'  = 2 + (p-1)n$
for any positive integer $n$.
Then there exist a cuspidal automorphic representation $\pi'$
of $\GL_2(\A_F)$, holomorphic of weight $(k',\ldots,k')$ such that
\begin{itemize}
\item if $v|p$, then $\pi'_v$ is unramified principal series;
\item $\ol{\rho}_{\pi'} \cong \rho$.
\end{itemize}
Suppose further that $\pi$ has prime-to-$p$ conductor dividing $\gn \subset \CO_F$,
and that $\psi$ is a finite order Hecke character of $\A_F^\times$ of conductor
dividing $\gn$, totally of parity $k'$, and satisfying $\det\rho = \overline{\psi\epsilon}$.
Then if $\rho$ is not badly
dihedral, we can choose $\pi'$ as above with conductor dividing $\gn$ and central
character $\psi^{-1}|\ |^{-k'}$.
\end{theorem}

\begin{remark}  Note that we may take $k' = p + 1$ in the conclusion, but in the
case $p=2$, this precludes using the Teichm\"uller lift of $\ol{\epsilon}^{-1}\det\rho$
for $\psi$ as this requires $k'$ to be even.
\end{remark}

\begpf  By Lemmas~\ref{hasse:Fp} and~\ref{hasse:Fq}, we see that
$[S_{(p-1,\ldots,p-1)}] \le [S_{(n(p-1),\ldots,n(p-1)}^{\otimes e}]$
for all $n,e \ge 1$, so by the arguments of Section~\ref{sec:lifting},
it suffices to prove that $\gm_{\rho\otimes\ol{\epsilon}^{-1}} \subset \T_{\CO_E}^S$ is in the
support of $S(U_1(\gn),V_{\{v|p\}})$ (with notation as in the proof
of Theorem~\ref{thm:variant1}, so in particular $V_{\{v|p\}}$ is the representation
of $U_1(\gn)\A_{F,f}^\times$ on which $U_p$ acts as 
$\otimes_{v|p} S_{v,(p-1,\ldots,p-1)}$, $U_x$ acts trivially for $x$
not dividing $p$, and $\A_{F,f}^\times$ acts via $\ol{\psi}$).

For $v|p$, define the representation $L_v$ of $\GL_2(\CO_{F,v})$ to be the
cokernel of the natural inclusion $\CO_E \to \Ind_{U_0(v)}^{\GL_2(\CO_{F,v})}\CO_E$.
We let $L_{\{v|p\}}$ denote the representation of $U_1(\gn)\A_{F,f}^\times$ on which $U_p$ acts
as $\bigotimes_{v|p} L_v$, $U_x$ acts trivially for $x$ not dividing $p$,
and $\A_{F,f}^\times$ acts via $\psi_\pi$, where $\psi_\pi^{-1}|\ |^{-2}$ is the central character
of $\pi$.  Note that $\ol{\psi}_\pi = \ol{\psi}$, so that
$L_{\{v|p\}} \otimes_{\CO_E} \Fpbar \cong V_{\{v|p\}}$;
moreover the induced inclusion
 $S(U_1(\gn),L_{\{v|p\}}) \otimes_{\CO_E} \Fpbar \to S(U_1(\gn),V_{\{v|p\}})$
is compatible with the natural action of $\T_{\CO_E}^S$.  Therefore it suffices to
prove that $\gm_{\rho\otimes\ol{\epsilon}^{-1}}$ is in the support of 
$S_0(U_1(\gn),L_{\{v|p\}}) = S(U_1(\gn),L_{\{v|p\}})/S^{\mathrm{triv}}(U_1(\gn),L_{\{v|p\}})$,
which in turn follows from it being in the support of 
$$S_0(U_1(\gn),L_{\{v|p\}})\otimes_{\CO_E}\C \cong 
  \bigoplus_{\Pi} (\Pi_f\otimes_{\CO_E} L_{\{v|p\}})^{U_1(\gn)\A_{F,f}^\times},$$
where the sum is over all cuspidal automorphic representations $\Pi = \Pi_f \otimes \Pi_\infty$
of $\GL_2(\A_F)$ such that $\Pi_\tau \cong D_{2,0}$  for all $\tau \in \Sigma$.
For $v|p$, $(\Pi_v \otimes_{\CO_E} L_v)^{\GL_2(\CO_{F,v})} \neq 0$ if and only
if $\Pi_v$ has unramified central character and conductor dividing $v$, so that
$(\Pi_f\otimes_{\CO_E} L_{\{v|p\}})^{U_1(\gn)\A_{F,f}^\times} \neq 0$
if and only if $\Pi$ has central character $\psi_\pi^{-1}$ and conductor dividing
$\gn \prod_{v|p} v$.  Furthermore note that $\gm_{\rho\otimes\ol{\epsilon}^{-1}}$ is in the support of
such a summand if and only if $\ol{\rho}_\Pi \cong \rho\otimes\ol{\epsilon}^{-1}$.
Finally our hypotheses ensure that $\Pi = \pi \otimes |\det|$ is exactly such an
automorphic representation.
\epf

\begin{remark}
Under the assumption that $\pi_v$ is an unramified principal series for all $v|p$ and other technical conditions ($p$ unramified in $F$, $[F:\Q]$ even or $F=\Q$), Theorem \ref{thm:variant2} is proved in Sections 2 and 4 of \cite{ek}.
\end{remark}

\end{document}